\documentclass{amsart}
\usepackage{latexsym}
\usepackage{graphicx}
\usepackage{epsfig}
\usepackage[active]{srcltx}
\usepackage{color}
\usepackage{amsmath,amssymb,amsfonts,mathrsfs, eufrak}%, txfonts}
\usepackage{pdfsync}
\usepackage[normalem]{ulem}

\def\dis{\displaystyle}

\begin{document}

% Theorems Lemmas ...
\newtheorem{thm}{Theorem}[section]
\newtheorem{lem}[thm]{Lemma}
\newtheorem{cor}[thm]{Corollary}
\newtheorem{df}[thm]{Definition}
\newtheorem{rem}[thm]{Remark}

% Proofs
\def\proof{{\noindent{\sc \bf Proof.\ }}}
\newcommand{\QED}{\hspace{1ex}\hfill$\Box$\vspace{2ex}}

% Greek letters, ...
\newcommand{\eps}{\varepsilon}
\newcommand{\al}{\alpha}
\def\t{\mathbb{\vartheta^0}}
\def\Ro{\mathbb{R}}

\newcommand{\Cbb}{\mathbb{C}}
\newcommand{\Dbb}{\mathbb{D}}
\newcommand{\Ebb}{\mathbb{E}}
\newcommand{\Ibb}{\mathbb{I}}
\newcommand{\Pbb}{\mathbb{P}}

% Domains
\newcommand{\Ome}{\mbox{$\Omega_\eps$}}
\newcommand{\ome}{\mbox{$\omega_\eps$}}
\newcommand{\Omte}{\mbox{$\widetilde\Omega_\eps$}}
\newcommand{\omte}{\mbox{$\tilde\omega_\eps$}}
\newcommand{\Bs}{{{B_\eps(s)}}}
\newcommand{\Peij}{{P^{\eps^\al}_{(i,j)}}}
\newcommand{\Seij}{{S^{\eps^\al}_{(i,j)}}}
\newcommand{\Steij}{{\tilde S^{\eps^\al}_{(i,j)}}}
\newcommand{\rij}{{r^{(i,j)}}}
\newcommand{\aij}{{a^{(i,j)}}}
\newcommand{\roij}{{\rho^{(i,j)}}}

\newcommand{\omin}{{\stackrel{\circ}{\omega}}}
\newcommand{\omine}{{\stackrel{\circ}{\omega}_\eps}}
\newcommand{\ominei}{{{\mathop{\omega}^{\,\,\,\circ\,\,i}}_{\!\!\eps}}}

\def\minint{\mathop{\int\hspace{-2.1ex}{\vspace{-0.5ex}-}}}

\newcommand\Haus{\mathcal{H}}
\newcommand\Leb{\mathcal{L}}

\newfont{\tenbss}{bbmss10}

\newcommand{\bsc}{\mbox{\tenbss c}}
\newcommand{\bsa}{\mbox{\tenbss a}}
\newcommand{\bss}{\mathfrak{s}}

\def\mint{{-}\hspace{-2.5ex}\int}
\def\fint{{-}\hspace{-2.5ex}\int}
\def\minint{\mathop{\int\hspace{-2.1ex}{\vspace{-0.5ex}-}}}
\newcommand{\wt}{\widetilde}
\newcommand{\weak}{\rightharpoonup}
\newcommand{\twoscale}{\rightharpoonup \hspace{-.25cm} \rightharpoonup}
\newcommand{\twoscalestrong}{\overset{\hspace{-.1cm}s}{\twoscale}}

%ENVIRONMENTS -------------------------------------------------------------
\newenvironment{remark}{\medskip \begin{rem}
\noindent \rm}{
\end{rem} %\begin{flushright} $\Box$n\end{flushright}
\medskip \par}

\title[Homogenization and dimension reduction for reinforced materials]{Non-local effects by homogenization\\ or 3D-1D dimension reduction \\
in elastic materials reinforced by stiff fibers}

\author[R. Paroni]{Roberto Paroni$^\ast$}
\thanks{$^\ast$
DADU, Universit\`{a} di
Sassari, Palazzo del Pou Salit, Piazza Duomo 6,
07041 Alghero, Italy, email: {\tt paroni@uniss.it} }
\author[A. Sili]{Ali Sili$^\dagger$}
\thanks{$^\dagger$D\'epartement de math\'ematiques, Universit\'e de Toulon, BP 20132, 83957, La Garde cedex, France        
                    \&
                     Institut de Math\'ematiques de Marseille (I2M, UMR 7373)
Centre de Math\'ematiques et Informatique (CMI), 
39 rue Joliot-Curie
13453 MARSEILLE Cedex 13
email: {\tt sili@univ-tln.fr}
}

\maketitle

%{\bf Version of \today.}

\noindent
{\bf \color{red} }

\medskip

\begin{abstract} 
We first consider an elastic thin heterogeneous cylinder of radius of order $\eps$: the interior of the cylinder is  occupied by a stiff material (fiber) that is surrounded  by a soft material (matrix). By assuming that the elasticity tensor of the fiber does not scale with $\eps$ and that of the matrix scales with $\eps^2$, we prove that the one dimensional model is a nonlocal system.

We then consider a reference configuration domain filled out by periodically distributed rods similar to those described above.  We prove that the homogenized model is a second order nonlocal problem.

In particular, we show that the homogenization problem is directly connected to the 3D-1D dimensional reduction problem. 
\end{abstract}

\medskip

{\small

\medskip

}
\section{\bf Introduction}\label{introduction}

Quite often nature combines two, or more, materials with complementary properties to generate a structural material whose performance and functionality  supersede those of monolithic materials,
\cite{WA}. Natural materials usually present a complex hierarchical structure, with characteristic dimensions spanning from the nanoscale to the macroscale, that is quite difficult to replicate in fabricated materials, \cite{WBS}. 
Nowadays, engineers and material scientists design and fabricate novel materials, inspired by natural materials, with a complex, but not hierarchical, structure. Two materials with completely different properties are appropriately arranged to create advanced functional materials that are lightweight, with remarkable strength, stiffness, and toughness, \cite{BM, PM, Ri}. The material properties of a two component material 
can be tuned by appropriately choosing the individual components, their morphology, their size, and arrangement, \cite{BM}. The understanding of the mechanical response of these advanced functional materials is an issue of paramount importance for industrial applications.

In the present paper, we study structures made up by a stiff and a soft material; the latter usually provides energy dissipation, toughness, ductility, and makes the structure lighter.  This combination is quite common in engineering and biological composites. 

We first derive, within the framework of linear elasticity, the elastic problem governing the motion of a rod composed of a stiff and a soft material. This is achieved by considering a sequence of problems posed on cylindrical reference configurations of diameters
that scale with a parameter $\eps$, and by taking a variational limit as $\eps$ approaches zero.   
The interior of the cylindrical regions is occupied by the stiff material (fiber), while the surrounding part by a soft material (matrix), whose modelling is achieved by scaling the elasticity tensor by $\eps^2$. In this way  the ratio between the components of the elasticity tensor of the fiber and those of the matrix  is $1/\eps^2$.
No assumption on the material symmetry of the body is made and the material could be inhomogeneous also within each of the stiff and the soft regions.

We show that the limit problem may be written equivalently as two independent systems: the first posed on the matrix, while the second takes into account  the elastic energy of only the fiber and the loads applied to the entire body. Mathematically speaking, the problem posed on the matrix is local, while the problem posed on the fiber is non-local since also the loads applied outside of the fiber region enter into the problem.
For any distribution of loads the problem is still non-local, since in order to get the global information about the limit displacement the problems posed on the matrix and on the fibers have to be both solved.
Leaving technicalities aside, the problem on the fiber determines a Bernoulli-Navier type of displacement, while the problem on the matrix determines the deviation from the Bernoulli-Navier type of displacement that takes place within the matrix region. 
The need of a displacement correction within the matrix region could be explained by the fact that the matrix is much more deformable than the fiber.

As a first step, to find the variational limit, we deduced a priori bounds on the displacements. These are easily obtained within the fiber region by means of Korn inequality, while a more intricate argument is needed within the matrix region. Indeed, since there is a loss of uniform ellipticity within the matrix region, only through a judicious application of several Poincar\'e and Korn inequalities, that also exploit the a priori bounds obtained within the fiber region, we are able to prove a priori bounds within the matrix region.

We also address a related homogenization problem. We consider a reference configuration domain $\Omega$ filled out by periodically distributed rods similar to those considered in the rod problem mentioned above. To completely fill  $\Omega$, we take rods with a square base of size $\eps$.
The body is therefore made of two regions: one occupied by a stiff material (fibers), while the second is filled by a soft material (matrix);  the ratio between the two elasticity components is still  $1/\eps^2$ as in the rod problem. We pass to the limit by using two-scale convergence,
this introduces a slow  variable $x$ and a fast  variable $y$. Both variables appear in the limit problem: the $x$ variable is used to describe the macroscopic behaviour while the $y$ variable the local, or microscopic, behaviour near each fiber. The limit problem could be recast in terms of only the macroscopic variable $x$, by taking advantage of the linearity of the problem; we refrain of doing it for the sake of brevity, see \cite{CS}. 
One of the main features of the limit problem is that locally each fiber behaves like the rod problem previously considered. As a consequence the homogenized problem is a second order problem: we recover a second-order material as a variational limit of a first order material, see also \cite{BB2,PS}. 
We finally emphasize that the techniques used connect the homogenization problem directly to the 3D-1D dimensional reduction problem. In particular, the non-locality already found in the rod problem still appears in the homogenization problem, for exactly the same reasons.

%%%%%%%%%%%%%%%%%%%%%%%%%%%%%

The literature on homogenization and on dimension reduction for the  diffusion or elasticity problems is huge, see for instance 
\cite{Allaire, BT,CL,CLM,EAG,Jar,FMP2008,Pan, Pa, AS1,AS2}.
We here describe only the works that are more related to ours.

In \cite{BB2} the homogenization of a periodic medium composed by stiff fibers surrounded by a soft matrix, similar to that considered in this paper, is studied. Both the fiber and the matrix are taken to be homogeneous and isotropic. Also, the elasticity tensor of the matrix scales with $\eps^2$, like in our case, but the scaling of the elasticity tensor of the fibers is proportional to $1/\eps^2$. The limit problem is non-local but it
does not depend on the extensional energy of the fibers.

The same homogenization problem with the same scaling as in \cite{BB2}, but in a more general context including  anisotropic and inhomogeneous materials, was considered in \cite{CS}. In addition to the non-local phenomenon arising in the limit and already proved in \cite{BB2}, it was shown in \cite{CS} that the anisotropy of the fibers occurs in the limit problem in the same way as in \cite{MS1} and \cite{MS2} leading to nonstandard terms. 
The phenomenon due to anisotropy occurs also in the present setting but to show it we would have to rewrite the problem in a reduced form. For the sake of brevity we do not give such formulation.

 We have chosen to scale the elasticity tensor of the matrix by $\eps^2$ and to keep the elasticity tensor of the fibers fixed,  i.e., it does not scale with $\eps$, among all the possible scalings, in order to have a non vanishing extensional energy of the fibers in the limit problem, differently from what happens in \cite{BB2,CS}. 

Also in \cite{BG} a similar homogenization problem, still for isotropic and homogeneous matrix and fibers, has been considered. The diameter of the fibers is allowed to scale with  size different from $\eps$, the elasticity tensor of the matrix is kept fixed, while that of the fiber may take several scalings. By taking the scaling of the latter equal to $1/\eps^2$ and by multiplying the full problem by $\eps^2$  our scaling is recovered but only with a load that scales at least as $\eps^2$: that is, a sequence of loads equivalent to a zero load.  

In both papers \cite{BB2,BG}, the fibers and the matrix are taken to be homogeneous and isotropic, while we do not make this assumption in the present paper.
It is well known, see \cite{CS, MS2, AS1,AS2}  that when the material under consideration presents high heterogeneities the nature of the limit problem is quite involved and a particular nonlocal phenomenon may appear in the limit. 

In the next section we fix the  notation and we formulate the problems. In Subsection  \ref{sec_rod} we define the rod problem while in Subsection \ref{sec_hom} we present the homogenization problem.
In Section \ref{results} we state our main results and we make some remarks. Section \ref{rod_proofs} and Section \ref{homo_proofs} are devoted to the proofs of the main results for the rod and homogenization problems, respectively.

\section{\bf Notation and setting of the problems}\label{problem}

As usual, and unless otherwise stated, Greek indices $ \alpha, \beta, \gamma,..., $  take values in the set $ \{1, 2\} $ while Latin indices $i, j , k, ...,$ in the set $ \{1, 2, 3\} $. 
Also, the summation convention for repeated indices is adopted throughout the paper.  

Often, for convenience, a $3\times3$ symmetric matrix $E$ is written as
\begin{equation}\label{tensorcomponents}
\begin{array}{cc}
E = \left(\begin{array}{cc}
 E_{\alpha\beta} & E_{\alpha 3}\\
E_{\alpha 3} & E_{3 3} 
\end{array}\right),
\end{array}
\end{equation}
where, with an abuse of notation, we denote by $E_{\alpha\beta}$ the 
$2\times2$ matrix whose components are $E_{\alpha\beta}$, and
similarly $E_{\alpha 3}$ denotes the  vector whose two components
are $E_{\alpha 3}$.
The product of two such matrices $E, K$ is defined by
\begin{equation}\label{tensorproduct}
E\cdot K :=   E_{ij} K_{ij} =:
 \left(\begin{array}{cc}
E_{\alpha\beta} & E_{\alpha 3}\\
E_{\alpha 3} & E_{3 3}
\end{array}\right)
\cdot
\left(\begin{array}{cc}
K_{\alpha\beta} & K_{\alpha 3}\\
K_{\alpha 3} & K_{3 3}
\end{array}\right).
\end{equation}

A generic point $ x \in \mathbb R^3$ is denoted by $ x = (x',  x_3)$, where $x'  = (x_1, x_2)\in\mathbb R^2$. 
The notation $x^R$ indicates the point obtained from $(x_1, x_2)$ by a rotation of  ${\pi}/{2}$ around the $x_3$-axis, i.e., $x^R := (- x_2, x_1)$. 
The gradient with respect to $x$ and $x'$ are denoted by $\nabla$ and $\nabla^{\prime}$, respectively, and the strain of a differentiable function $u$ is denoted by 
$$Eu:={\rm sym}\nabla u:=\frac 12 (\nabla u+\nabla u^T).$$
Throughout the paper $\partial_i$ denotes the partial derivative with respect to $x_i$.

In the homogenization problem we shall have, besides the variable $x$, also a variable $y$. We denote by 
$$\partial^y_\alpha\varphi :=\frac{\partial \varphi}{\partial {y_\alpha}},
$$
for any differentiable function $\varphi$.
Also,
\begin{equation}\label{Eyx}
(E^y\varphi)_{\alpha \beta}:=\frac 12 (\partial^y_\alpha\varphi_\beta+\partial^y_\beta\varphi_\alpha),
\qquad
(E^{yx}\varphi)_{\alpha3}:=\frac 12 (\partial^y_\alpha \varphi_3+\partial_3 \varphi_\alpha),
\end{equation}
for every differentiable vectorial field $\varphi$.

\subsection{The rod problem}\label{sec_rod}

Let $\omega, \omin \subset \Ro^2$ be two open, bounded, simply connected sets with Lipschitz boundaries, $\partial\omega$ and $\partial \omin$, respectively, such that $\omin\subset\subset \omega$. \par
We assume the $x_\alpha$-axes to be centered in the center of mass of $\omin$ so to have that
\begin{equation}\label{centermass}
\int_\omin x_\alpha  \ dx' = 0.
\end{equation}

Let $\ell>0$, $I:=(0,\ell)$, and
$$\Omega:=\omega\times I,  \quad F:=\omin\times I,\quad M:=\Omega\setminus \bar F,$$
where $\bar F$ denotes the closure of $F$.
Physically, we can think of $\Omega$ as the reference configuration of the rod, $F$ as the region occupied by the fiber (strong material), and $M$ as the region occupied by the matrix (soft material). If $x\in \Omega$ we have that $x'\in \omega$ and $x_3 \in I$.

For every $\eps>0$, let $R^\eps\in \Ro^{3\times 3}$ be the diagonal matrix whose entries are $\eps, \eps,$ and $1$, i.e., $R^\eps=\mbox{diag}(\eps,\eps,1)$.
The scaled gradient $H^\eps u$ and the scaled strain $E^\eps u$ are defined by
$$H^\eps u:= (R^\eps)^{-1}\nabla u (R^\eps)^{-1}, \quad E^\eps u:= (R^\eps)^{-1}E u (R^\eps)^{-1}=
{\rm sym}  H^\eps u.$$
In components we have
\begin{equation}\label{components}
\begin{array}{cc}
\displaystyle (H^\eps u)_{\alpha\beta}=\frac 1{\eps^2} \partial_\beta u_{\alpha}, &
\displaystyle (H^\eps u)_{\alpha3}=\frac 1{\eps} \partial_3 u_{\alpha}, \\[2ex]
\displaystyle (H^\eps u)_{3\beta}=\frac 1{\eps} \partial_\beta u_3, &
\displaystyle (H^\eps u)_{33}= \partial_3 u_{3}.
\end{array}
\end{equation} \par

Let $\Cbb$ be a fourth-order symmetric tensor field defined on $\Omega$. We assume that $\Cbb$ fulfills the following assumptions: for all $i,j,k,l$, for all $3\times 3$ symmetric matrices  $E$, and for almost all $x \in \Omega$

\begin{equation}\label{tensor}
\left\{
\begin{array}{ll}
\displaystyle a) \  \Cbb_{ijklÊÊ}(x) = \Cbb_{jikl}(x) = \Cbb_{klij}(x),\\
 \displaystyle b)  \ \Cbb_{ijklÊÊ} \in L^\infty(\Omega), \\
c) \ \hbox{there exists} \ m > 0, \, \, \Cbb_{ijklÊÊ}(x) E_{kl }E_{ij} \geq m \  E_{ij }E_{ij}.
\end{array}
\right.
\end{equation}

We assume the rod to be clamped at both ends of the cylinder $\Omega$, we thus set
 $$H^1_{dd}( \Omega):=\left\{u\in H^1( \Omega; \Ro^3):
u(x',0)=u(x',\ell)={0} \right\},$$ 
and subjected to body forces $f\in L^2(\Omega;\Ro^3)$.
In the notation above the two $d$'s in $H^1_{dd}$ are used to recall that on both ends of the beam we impose Dirichlet boundary conditions.

By $u^\eps\in H^1_{dd}( \Omega)$ we denote the solution of
\begin{equation}\label{mainpb}
\left\{
\begin{array}{ll}
\displaystyle u^\eps\in H^1_{dd}( \Omega)\\
\displaystyle  \int_\Omega (\chi_F+\eps^2\chi_M)\Cbb(x) E^\eps u^\eps\cdot E^\eps \varphi\,dx=
\int_\Omega f(x)\cdot \varphi(x)\,dx\quad \forall \varphi \in H^1_{dd}( \Omega).
\end{array}
\right.
\end{equation}
By the previous assumptions on the tensor $\Cbb$ and the body force $f$, problem \eqref{mainpb} is well-posed by the Lax-Milgram Theorem.

Problem \eqref{mainpb} is a rescaled elasticity problem, as explained in Remark \ref{omegaeps} below, of a body whose elasticity tensor is equal to $\Cbb$ on $F$ and equal to $\eps^2 \Cbb$ on $M$. Thus, the elasticity of the fiber, i.e., the elasticity of the material occupying the region $F$, is $1/\eps^2$ larger than the elasticity of the matrix (material in $M$). We point out that our assumptions allow both the fiber and the matrix to be inhomogeneous and to be fully anisotropic.

One of the aims of the present paper is to deduce the variational limit of \eqref{mainpb} as $\eps$ goes to zero.
\begin{rem}\label{omegaeps}
{\rm
 Problem \eqref{mainpb} is the variational formulation in the fixed domain $\Omega$ of an elasticity
 problem posed in the variable thin domain $ \Omega_\eps$, see \cite{CiDe}. Indeed, let 
$$\Omega_\eps:=\eps\omega\times I,  \quad F_\eps:=\eps\omin\times I,\quad M_\eps:=\Omega_\eps\setminus \overline F_\eps.$$

Let $r^\eps:\Omega\to\Omega_\eps$ be defined by $r^\eps(x):=(\eps x_1,\eps x_2,x_3)$. Then
$R^\eps = \nabla r^\eps$.
For any $v:\Omega\to \Ro^3$ we define $\hat v:\Omega_\eps\to \Ro^3$ by
$$
\hat v:= (R^\eps)^{-1} v\circ ( r^\eps)^{-1},
$$
in components
$$
\hat v_\alpha(\eps x_1,\eps x_2,x_3)=\frac 1\eps v_\alpha(x), \qquad
\hat v_3(\eps x_1,\eps x_2,x_3)= v_3(x).
$$
We note that
$$\nabla \hat v=(R^\eps)^{-1} \nabla v\circ (r^\eps)^{-1}(R^\eps)^{-1} = H^\eps v\circ (r^\eps)^{-1},\quad
E\hat v = E^\eps v\circ (r^\eps)^{-1},
$$
With this notation,  problem \eqref{mainpb} rewrites as
\begin{equation}\label{pbvariabledomain}
\left\{
\begin{array}{ll}
\displaystyle \hat u^\eps\in H^1_{dd}( \Omega_\eps)\\
\displaystyle  \int_{\Omega_\eps} (\chi_{F_\eps}+\eps^2\chi_{M_\eps})\Cbb^\eps(x) E \hat u^\eps\cdot E\hat \varphi\,dx=
\int_{\Omega_\eps} f^\eps(x)\cdot \hat\varphi(x)\,dx, \quad
 \forall \hat\varphi\in H^1_{dd}( \Omega_\eps),
\end{array}
\right.
\end{equation}
where
$$
\Cbb^\eps:=\Cbb\circ (r^\eps)^{-1},\qquad f^\eps:=R^\eps f\circ (r^\eps)^{-1}.
$$
From the previous definition, the components of the load have two different
scalings:
\begin{equation}\label{loadscaling}
f^\eps_\alpha=\eps f_\alpha\circ (r^\eps)^{-1},\quad f^\eps_3= f_3\circ (r^\eps)^{-1}.
\end{equation}
}
\end{rem}

\subsection{The homogenization problem}\label{sec_hom}

The reference configuration of the elastic body is now assumed to be the parallelepiped  
 $ \Omega = \omega \times (0, \ell) = \omega \times I  $, where $\ell>0$ and $\omega$ is the square  $ \omega = (-\ell, \ell)^2$.
 We consider a sequence $\eps$ approaching zero such that
 $$\omega = \bigcup_{ i \in I_\eps} \omega_\eps^i, \quad \mbox{where} \, \,  \omega_\eps^i =  \eps Y + \eps (i_1, i_2) \, \, \mbox{and} \, \, Y = ( - \frac12, \frac12) ^2.$$ 
 Let $ \Omega_\eps^i$ and $I_\eps$ 
 be defined by
  $$ \Omega_\eps^i :=  \omega_\eps^i \times (0, \ell)$$
  and
 $$ I_\eps := \{ i = (i_1, i_2) \in \mathbb{Z}^2, \, \,  \omega_\eps^i \subset \omega \}. $$
 Let $D$ be the disk centered at the origin with radius  $ 0 < r < \frac12$.  Set 
 $$D_\eps^i:= \eps D + \eps (i_1, i_2).$$ 
 The set $F_\eps$ of fibers contained in $\Omega$ and its complementary set in $\Omega$,  the matrix $M_\eps$,  are defined by:    
   \begin{equation}\label{thegeometryhomogenization}
\left\{
\begin{array}{ll}
F_{\eps} := \displaystyle \bigcup_{i
\in I_{\eps} } F_{\eps}^i, \, \, \mbox{where} \,  F_{\eps}^i := D_\eps^i \times (0, \ell) \\
 \displaystyle  M_{\eps}^i := \Omega_{\eps}^i \setminus \bar F_{\eps}^i, \quad M_\eps := \Omega \setminus \bar F_\eps = \displaystyle \bigcup_{i
\in I_{\eps} } M_{\eps}^i, \\
\displaystyle \Omega = \bigcup_{{i \in I_{\eps}}}
\Omega_{\eps}^i = F_\eps \cup M_\eps.
\end{array}
\right.
\end{equation}
   
  The microscopic variable will be denoted by $ y=(y_1,y_2)$ (notice that there is no fast variable in the vertical direction) and, as before, $y^R= (-y_2, y_1)$. 
  
 We now assume that the elasticity tensor $\Cbb$ is such that: 
   \begin{equation}\label{tensorhomogenization}
\left\{
\begin{array}{ll}
\displaystyle a) \  \Cbb_{ijklÊÊ}(x,y) = \Cbb_{jikl}(x,y) = \Cbb_{klij}(x,y),\\
 \displaystyle b)  \ \Cbb_{ijkl } \in L^\infty (\Omega; C_{\#}(Y)), \\
c) \ \mbox{there exists} \ m > 0, \, \, \Cbb_{ijklÊÊ}(x,y) E_{kl }E_{ij} \geq m \  E_{ij }E_{ij},
\end{array}
\right.
\end{equation}
for almost all $(x,y) \in \Omega \times Y$, for
all $i,j,k,l$, and for all symmetric $3 \times 3$ matrices $E$, and 
 where $ C_{\#}(Y)$ denotes the space of continuous $Y$-periodic functions on $\Ro^2$. 
 
As for the rod problem, we assume the body to be made up from a set of strong fibers surrounded by a soft matrix, with ratio between the elasticity tensor of the soft part and the strong part being equal to $1/{\eps^2}$. Therefore,
  the variational problem modeling the equilibrium of the heterogeneous elastic medium at the microscopic level may be written as:
  \begin{equation}\label{mainpbhomogenization}
\left\{
\begin{array}{ll}
\displaystyle u^\eps \in H^1_{dd}( \Omega),\\
\displaystyle  \int_\Omega (\chi_{F_\eps}+\eps^2\chi_{M_\eps})\Cbb(x, \frac {x'}{ \eps}) E u^\eps\cdot E \varphi\,dx =
\int_\Omega \eps f_\alpha(x, \frac{x'}{\eps})\cdot \varphi_\alpha(x)  \\
\displaystyle\hspace{7cm}+ f_3(x, \frac{x'}{\eps})\cdot \varphi_3(x) \,dx, \\
\displaystyle   \forall \varphi \in H^1_{dd}( \Omega).
\end{array}
\right.
\end{equation}
 
Note that compared to \eqref{pbvariabledomain}, only the shape of the domain $\Omega$ has changed. Note also that we  have assumed, in \eqref{mainpbhomogenization}, the loading $f(x,y)$ to depend on both the macroscopic variable $x$ and the microscopic variable $y$, and that
the scaling of the loads coincides with that of  \eqref{pbvariabledomain}, see \eqref{loadscaling}.

By assuming $f_i \in L^2 (\Omega; C_{\#}(Y))$, we deduce the existence and the uniqueness of the solution $u^\eps$,  for each fixed $ \eps > 0$, by means of the Lax-Milgram Theorem.  

To derive the limit problem, we will use compactness properties related to the two-scale convergence, see \cite{Allaire} and \cite{Nguetseng}. Recall that a sequence $u^\eps \in L^2(\Omega)$ (weakly) two-scale converges to a function $u \in L^2(\Omega \times Y)$, written $u^\eps\twoscale u$, if
 $$ \dis \int_\Omega u^\eps(x) \varphi(x, \frac {x'}\eps) \ dx \to  \int_{\Omega \times Y} u(x,y)  \,  \varphi(x,y) \ dx dy \quad \forall \, \varphi \in L^2(\Omega;  C_{\#}(Y)). $$  

A sequence $u^\eps \in L^2(\Omega)$ is said to strongly two-scale converge to  $u \in L^2(\Omega \times Y)$, written $u^\eps\twoscalestrong u$, if
 $$ \dis \int_\Omega u^\eps(x) v^\eps(x) \ dx \to  \int_{\Omega \times Y} u(x,y)  \,  v(x,y) \ dx dy$$
 for any  $v^\eps \in L^2(\Omega)$ that weakly two-scale converges to  $v \in L^2(\Omega \times Y)$.
 
 We recall a few properties of the two-scale convergence that will be useful in the analysis below:
 
 \begin{itemize}
 \item any bounded sequence in $L^2(\Omega)$ admits a two-scale converging subsequence; 
 \item if $u^\eps$ is bounded in $H^1(\Omega)$ then there exist $u \in H^1(\Omega)$ and $u^1 \in L^2(\Omega; H^1_{\#} (Y))$ such that, up to a subsequence, 
  $$ u^\eps \twoscale u, \, \, \nabla u^\eps \twoscale \nabla_x u + \nabla_y u^1;$$
\item if $u^\eps$ is bounded in $L^2(\Omega)$ and if $ \eps \nabla u^\eps$ is bounded in $L^2(\Omega; \Ro^3)$, then there exists $u^1 \in L^2(\Omega; H^1_{\#} (Y))$ such that, up to a subsequence, 
  $$ u^\eps \twoscale u^1, \quad \eps \nabla u^\eps \twoscale  \nabla_y u^1.$$
 \item if  $u^\eps \in L^2(\Omega)$ weakly two-scale converges to  $u \in L^2(\Omega \times Y)$ and
$$
\lim_{\eps\to 0}\|u^\eps\|_{L^2(\Omega)}=\|u\|_{L^2(\Omega\times Y)}
$$
 then $u^\eps \twoscalestrong u$, see \cite{ZY}.
\end{itemize}

\section{Main results}\label{results}
The main results of the paper are stated in the present section. 

\subsection{The rod problem}

We first define a few spaces that will be useful in the results stated below.
The space of Bernoulli-Navier type of displacements satisfying the prescribed boundary conditions 
is denoted by
\begin{eqnarray*}
\mathcal{BN}_{dd}(F)&:= &\{u\in H^1_{dd}(F):(Eu)_{i\alpha}=0 \hbox{ a.e. in } F\}\\
& = &
\{u\in H^1_{dd}(F):\exists \xi_\alpha\in H^2_{0}(0,\ell),\exists\xi_3\in H^1_{0}(0,\ell) \mbox{ s.t. }\\
& & \hspace{3cm}u_\alpha=\xi_\alpha, \, u_3 = \xi_3-x_\alpha \xi_\alpha', 
 \hbox{ a.e. in }  F \}.
\end{eqnarray*}
We note that the elements of $\mathcal{BN}_{dd}(F)$ may be naturally extended to elements of $\mathcal{BN}_{dd}(\Omega)$.
The set of functions orthogonal to the bi-dimensional rigid displacements
is denoted by
$$
 \mathcal{RD}^\perp_{2}(F):=
\{(w_1,w_2,0): w_\alpha\in L^2(I;H^1(\omin)), \int_\omin x_1w_2-x_2w_1\,dx'=0
\}
$$
while the set of twisting and cross-sectional warping displacements is 
$$
\mathcal{R}_{dd}(F):=\{ (-x_2\vartheta, x_1\vartheta,v_3):\ \hbox{ with} \, \vartheta \in H^1_0(0,\ell) \, \hbox{and } \  v_3\in L^2(I;H^1_m(\omin))\}.
$$
where
$$
H^1_m(\omin):=\{\varphi\in H^1(\omin): \mint_\omin\varphi\,dx'=0\},
$$
For brevity, we set
$$
\mathcal{U}:=\mathcal{BN}_{dd}(F), \quad
\mathcal{V}:=\mathcal{R}_{dd}(F), \quad
\mathcal{W}:=\mathcal{RD}^\perp_{2}(F),
$$ 
and
$$
\mathcal{Z}:=\{z: z_i\in L^2(I;H^1(\omega)), z_i=0 \mbox{ a.e.\ in } F \}.
$$
The following theorem is our main result for the rod problem.

\begin{thm}\label{Mainresult}
There exists $(u,v,w,z) \in \mathcal{U} \times \mathcal{V}  \times \mathcal{W} \times \mathcal{Z}$ such that the sequence of solutions $u^\eps$ of problem \eqref{mainpb} fulfills the following convergences
$$
u^\eps_\alpha\to u_\alpha \mbox{ in }H^1(\Omega),$$
$$
u^\eps_3\to  u_3+z_3 \mbox{ in } L^2(I;H^1(\omega)),
$$
$$
E^\eps u^\eps\chi_F\to
\left(\begin{array}{cc}
(Ew)_{\alpha\beta} & (Ev)_{\alpha 3}\\
(Ev)_{\alpha 3} & (Eu)_{3 3}
\end{array}\right)\chi_F
  \quad  \mbox{in } \, L^2(\Omega),
$$ 
$$ E^\eps u^\eps\chi_M\to 
\left(\begin{array}{cc}
(Ez)_{\alpha\beta} & \frac 12 \partial_\alpha z_3\\
\frac 12 \partial_\alpha z_3 & 0
\end{array}\right)
\chi_M \quad   \mbox{in } \, L^2(\Omega).   $$ 
The limit $(u,v,w,z)$ is the unique solution of the problem
\begin{equation}\label{Thelimitpb}
\left\{
\begin{array}{ll}
\displaystyle (u,v,w,z)\in  \mathcal{U}\times\mathcal{V}\times\mathcal{W}\times\mathcal{Z}, \\[1ex]
\displaystyle  
\int_\Omega \Cbb(x)
\left(\begin{array}{cc}
(Ew)_{\alpha\beta} & (Ev)_{\alpha 3}\\
(Ev)_{\alpha 3} & (Eu)_{3 3}
\end{array}\right)
\cdot
\left(\begin{array}{cc}
(E\bar w)_{\alpha\beta} & (E\bar v)_{\alpha 3}\\
(E\bar v)_{\alpha 3} & (E\bar u)_{3 3}
\end{array}\right)
\chi_F\,dx
\, + \\[4ex]
\hspace{1cm}\displaystyle  
\int_\Omega \Cbb(x)
\left(\begin{array}{cc}
(Ez)_{\alpha\beta} & \frac 12 \partial_\alpha z_3\\
\frac 12 \partial_\alpha z_3 & 0
\end{array}\right)
\cdot
\left(\begin{array}{cc}
(E\bar z)_{\alpha\beta} & \frac 12 \partial_\alpha \bar z_3\\
\frac 12 \partial_\alpha \bar z_3 & 0
\end{array}\right)
\chi_M\,dx
=  \\[4ex] 
\hspace{4,5cm}\displaystyle  
\int_\Omega f_\alpha(x)\bar u_\alpha(x_3)+f_3(x)(\bar u_3(x)+\bar z_3(x))\,dx, \\[4ex] \forall(\bar u,\bar v,\bar w,\bar z)\in \mathcal{U}\times\mathcal{V}\times\mathcal{W}\times\mathcal{Z}.
\end{array}
\right.
\end{equation}
\end{thm}
\vspace{0.5cm}

\begin{rem}\label{onthemainresult} 
{\rm
Some remarks are in order.
\begin{itemize}
\item The space $ \mathcal{U}\times\mathcal{V}\times\mathcal{W}\times\mathcal{Z}  $ equipped with the norm 
\begin{eqnarray*} \displaystyle & \parallel (u,v,w,z) \parallel^2 =  \displaystyle \sum_{\alpha \beta}  
( \parallel (Eu)_{33} \parallel^2_{L^2(\Omega)}   + \parallel (Ev)_{\alpha3} \parallel^2_{L^2(\Omega)} + \parallel (Ew)_{\alpha \beta} \parallel^2_{L^2(\Omega)}   \\
 & + \parallel (Ez)_{\alpha \beta} \parallel^2_{L^2(\Omega)} + \parallel \partial_\alpha z_3 \parallel^2_{L^2(\Omega)}    )  
 \end{eqnarray*}  
 is a Hilbert space (for an analogous proof in the homogenization framework one can refer to  \cite{CS}, Lemma 4.5). Therefore, the well-posedness of the limit problem \eqref{Thelimitpb} is easily obtained from the Lax-Milgram Theorem since $ f \in (L^2(\Omega))^3$ and $\Cbb$ is strongly coercive, \eqref{tensor}. 

\item The limit displacement can be described, in words, as follows: it is a Bernoulli-Navier displacement with the axial component augmented over the matrix region. Said differently, $z_3$ takes into account the deviation from a Bernoulli-Navier type of displacement in the matrix region. 
The need of this correction over the matrix region could be explained by the fact that the matrix is much more deformable than the fiber, because of the different stiffnesses. 

\item The system posed on the fiber $F$ coincides with a problem obtained in \cite{MS1}; this is not surprising since the sequence $u^\eps$ has $(E^\eps u^\eps)_{ij}$ bounded in $L^2(F)$, as in \cite{MS1}. In particular, it was shown in \cite{MS2}, and also in \cite{CS} in the homogenization framework, that the displacements $(v,w)$ may be expressed in terms of the Bernoulli-Navier displacement $u$. For brevity, we will not reproduce here the proofs of these statements.

\item Choosing successively $(\bar u, \bar v, \bar w) = (0, 0, 0)$ and then $ \bar z =0$, we see that problem \eqref{Thelimitpb} may be written equivalently as two  systems: the first one posed on the matrix $M$ with unknown  $z$ and the second essentially posed on the fiber $F$ with unknown $(u,v,w)$. Indeed, the second problem considers  the elastic energy of the fiber and the loads applied to the fiber \emph{and} the matrix. Thanks to the Bernoulli-Navier structure of the displacement $u$ one can redefine the loads and rewrite the second problem over the domain $F$ only; the redefined loads will depend on the original loads that were acting also on the matrix.

The well-posedness of these two problems may be established independently, but
these problems present a non-local effect: the displacements in the fiber are influenced also by
the loads applied to the matrix.

Also,  to describe the limit displacements, even in the region $M$, the solution of both systems are needed, as, for instance, stated by the convergence $u^\eps_3\to  u_3 + z_3 \mbox{ in }L^2(I;H^1(\omega))$, note that $z_3=0$ in $F$, see Theorem \ref{Mainresult}.
\end{itemize}
}
\end{rem}

The non-local effect mentioned in the previous remark is further studied in the following theorem, which also gives the limit of the strains defined on the variable domain.
\begin{thm}\label{nonlocaleffect}
Let $(u,v,w,z)$ be the solution of \eqref{Thelimitpb}. Let
$$E_f :=  \left(\begin{array}{cc}
(Ew)_{\alpha\beta} & (Ev)_{\alpha 3}\\
(Ev)_{\alpha 3} & (Eu)_{3 3}
\end{array}\right),  \quad
 E_m := \left(\begin{array}{cc}
(Ez)_{\alpha\beta} & \frac 12 \partial_\alpha z_3\\
\frac 12 \partial_\alpha z_3 & 0
\end{array}\right).  $$
Then the sequence of solutions $\hat u^\eps$ of \eqref{pbvariabledomain} is such that
$$ \mint_{\eps \omin}E \hat u^\eps  \ dx'  \to  \mint_{ \omin} E_f  \ dx'  , \qquad  \mint_{\eps (\omega \setminus \omin) } E \hat u^\eps  \ dx'  \to  \mint_{ \omega \setminus \omin} E_m  \ dx' $$
 in $L^2(I;\Ro^{3 \times 3}),$
$$\mint_{\eps\omega}{\eps} \hat u^\eps_\alpha\ dx' \to  u_\alpha(x_3) \mbox{ in } \, L^2(I),  $$
$$
\mint_{\eps\omega} \hat u^\eps_3 \ dx' \to  U(x_3) := \mint_\omega (u_3 + z_3) \ dx' \mbox{ in } \, L^2(I).
$$
$U$ may be written as 
$$
U(x_3) = \mint_{\omega} u_3 \ dx' + m_0(x_3) \mint_{\omega \setminus \omin} f_3 \ dx' + m_{00}(x_3),  
$$ 
with $m_0$ and $m_{00}$ given by
$$ m_0 (x_3) = \frac{1}{\vert \omega \vert}\int_{\omega \setminus \omin} z^0_3 \ dx', \quad   m_{00} (x_3) = \frac{1}{\vert \omega \vert}\int_{\omega \setminus \omin} z^{00}_3 \ dx',$$ 
where $z^0$ and $z^{00  }$ are respectively the solutions of the following problems
\begin{equation}\label{elementarysolution}
 \left\{ 
\begin{array}{ll}
z^0 \in L^\infty(I; \mathcal{Z}^0)   \\[2ex]   
 \displaystyle \int_{\omega \setminus \omin} \Cbb(x)
\left(\begin{array}{cc}
(Ez^0)_{\alpha\beta} & \frac 12 \partial_\alpha z^0_3\\
\frac 12 \partial_\alpha z^0_3 & 0
\end{array}\right)
\cdot
\left(\begin{array}{cc}
(E\bar z)_{\alpha\beta} & \frac 12 \partial_\alpha \bar z_3\\
\frac 12 \partial_\alpha \bar z_3 & 0
\end{array}\right)
\ dx'
=   \\[4ex] \hspace{6,5cm}\displaystyle 
\int_{\omega \setminus \omin} \bar z_3(x'))\,dx' 
\quad\mbox{ a.e.\ in } I\\[3ex]  
 \forall \bar z \in \mathcal{Z}^0 := \{ z : z_i \in H^1(\omega \setminus \omin), \, z_i = 0 \, \mbox{on} \ \partial \omin \},
\end{array}
\right.
\end{equation}

\begin{equation}\label{secondsolutiononZ}
 \left\{ 
\begin{array}{ll}
z^{00} \in \mathcal{Z},  \\[1ex]   
 \displaystyle \int_{M} \Cbb(x)
\left(\begin{array}{cc}
(Ez^{00})_{\alpha\beta} & \frac 12 \partial_\alpha z^{00}_3\\
\frac 12 \partial_\alpha z^{00}_3 & 0
\end{array}\right)
\cdot
\left(\begin{array}{cc}
(E\bar z)_{\alpha\beta} & \frac 12 \partial_\alpha \bar z_3\\
\frac 12 \partial_\alpha \bar z_3 & 0
\end{array}\right)
\ dx
=   \\[4ex] \hspace{5cm}\displaystyle 
\int_{M} (f_3(x) - \mint_{\omega\setminus \omin} f_3(x) \ dx' )\bar z_3(x) \,dx \\ 
 \forall \bar z \in \mathcal{Z}.
\end{array}
\right.
\end{equation}
\end{thm}

\begin{rem}\label{onthenonlocaleffect}
\mbox{ }
{\rm
\begin{itemize}
\item
The last convergence in Theorem 3.3  emphasises the non-local effect. 
This can be easily seen in the particular case in which the loading term $f_3$  depends only on $x_3$. In this case, we have $z_3^{00}=0$ and hence
$$
\mint_{\eps\omega} \hat u^\eps_3 \ dx' \to  U(x_3) = \mint_\omega u_3  \ dxÕ+m_0(x_3) \mint_{\omega \setminus \omin} f_3 \ dx' \mbox{ in } \, L^2(I).
$$
Thus the axial displacement in the limit problem is the sum of two terms:
\begin{itemize}
\item $\displaystyle \mint_\omega u_3 \ dx' $ which comes by solving a problem over the fiber but by taking into account also the loads over the matrix;
\item $ \displaystyle m_0(x_3) \mint_{\omega \setminus \omin} f_3 \ dx'$ which is deduced by solving a displacement problem over the matrix region. Note that from \eqref{elementarysolution} it follows that $m_0(x_3)>0$ for a.e. $x_3\in I$.
\end{itemize}
 \item The displacement 
 $$
U(x_3) = \mint_{\omega} u_3 \ dx' + m_0(x_3) \mint_{\omega \setminus \omin} f_3 \ dx' + m_{00}(x_3),  
$$ 
 has been split in the sum of three terms: the first is the average of the Bernoulli-Navier displacement, the second takes into account the correction due to  uniform loads over the matrix region and the last the presence of variable loads with zero average.
\item All the previous results remain true when the cylinder $\Omega^\eps  $ is assumed to be clamped on only one of its ends; the limit system is the same as \eqref{Thelimitpb} with the spaces 
$\mathcal{U}$ and $\mathcal{V}$ defined as before but with the space $H^1_{dd}(\Omega)$ 
replaced by 
$$  H^1_{d}( \Omega):=\left\{u\in H^1( \Omega;\Ro^3):
u(x',\ell)= 0 \right\},$$
while the spaces $\mathcal{W}$ and $\mathcal{Z}$ remain unchanged, since the boundary conditions do not enter in their definition.
\end{itemize}
} 

\end{rem}

\subsection{The homogenization problem}

 The spaces we will use to describe the homogenized problem are the following (below, $h$ stands for homogenization):

\begin{eqnarray}
\mathcal{U}^h \!\!\!\!\!&:=&\!\!\!\!\! \{(u_i,u^1_\alpha, \vartheta):  u_\alpha \in H^1(\Omega)\cap L^2(\omega;H_0^2(I)), \vartheta\in L^2(\omega;H^1_0(I)),\nonumber\\
&&\hspace{1,7cm}   u^1_\alpha \in  L^2(\Omega; H^1_{\#}(Y))\cap L^2(\Omega; H^1_{m}(D))\cap L^2(\omega\times D;H^1_0(I)), \nonumber\\
&&\hspace{1,7cm} \mbox{and }  u_3\in  L^2(\Omega; H^1_{\#}(Y)),  u_3(\cdot,y)\in L^2(\omega;H^1_0(I)) \mbox{ for a.e. } y\in Y, \nonumber\\
& & \hspace{1,7cm} \mbox{such that }\label{defUh}\\
& &\hspace{1,7cm} \dis u_3 (x,y) = - y_\alpha \partial_3u_\alpha(x) +\mint_D u_3 (x,y)\,dy, \, \mbox{ for a.e. } (x,y)\in \Omega \times D\nonumber\\
& & \hspace{1,7cm}u^1_\alpha(x,y)=- y_\beta \partial_\beta u_\alpha (x) + \vartheta(x) y^R_\alpha%
\mbox{ for a.e.}\ (x,y)\in\Omega\times D\nonumber\},
\end{eqnarray}

$$
\mathcal{V}^h:= L^2(\Omega; H^1_m(D)),
$$
where
$$
H^1_m(D):=\{\varphi\in H^1(D): \mint_D\varphi\,dy=0\},
$$
and
$$
 \mathcal{W}^h:=
\{(w_1,w_2,0): w_\alpha\in L^2(\Omega ; H^1(D)), \int_D y_1 w_2 - y_2 w_1\,dy=0, \, \, \mbox{a.e.} \  x \in \Omega
\}.
$$

The following theorem, analogous to Theorem \ref{Mainresult}, is our main result for the homogenization problem.

\begin{thm}\label{Mainresulthom}
There exist
$(u_i,u^1_\alpha,\vartheta) \in \mathcal{U}^h$,  $v_3 \in \mathcal{V}^h$, and $w \in \mathcal{W}^h$
such that the sequence of solutions $u^\eps$ of problem \eqref{mainpbhomogenization} fulfills the following convergences
$$
\eps u^\eps_\alpha \weak u_\alpha, \mbox{ in }H^1(\Omega),
\quad
\eps \nabla^\prime u^\eps_\alpha \twoscale \nabla^\prime u_\alpha +\nabla^\prime_y u^1_\alpha, 
$$
$$
u^\eps_3 \twoscale  u_3,
$$
Moreover, setting $v_\alpha:=y_\alpha^R \vartheta$, we have
$$
(E u^\eps)\chi_{F_\eps} \twoscalestrong 
\left(
\begin{array}{cc}
(E^y w)_{\alpha\beta} & (E^{yx}v)_{\alpha3}\\
(E^{yx}v)_{\alpha3} & (Eu)_{33}
\end{array}
\right)
\chi_D(y),
$$ 
$$ \eps(E u^\eps)\chi_{M_\eps} \twoscalestrong  
\left(
\begin{array}{cc}
(Eu + E^y u^1)_{\alpha\beta} & (E^{yx}u)_{\alpha3}\\
(E^{yx}u)_{\alpha3} & 0
\end{array}
\right)
\chi_{Y \setminus D} (y).   
$$ 
The limit $(u_i,u^1_\alpha,\vartheta,v_3,w)$ is the unique solution of the problem
\begin{equation}\label{Thelimitpbhom}
\left\{
\begin{array}{l}
\displaystyle (u_i,u^1_\alpha,\vartheta,v_3,w)\in  \mathcal{U}^h\times\mathcal{V}^h\times\mathcal{W}^h, 
\\[1ex]
\displaystyle  
\int_\Omega\int_D 
\Cbb
\left(
\begin{array}{cc}
(E^y w)_{\alpha\beta} & (E^{yx}v)_{\alpha3}\\
(E^{yx}v)_{\alpha3} & (Eu)_{33}
\end{array}
\right)
\cdot
\left(
\begin{array}{cc}
(E^y\bar w)_{\alpha\beta} & (E^{yx}\bar v)_{\alpha3}\\
\displaystyle  
(E^{yx}\bar v)_{\alpha3} & \partial_3\bar u_3
\end{array}
\right)
\,dydx
\\[15pt]
\displaystyle  
 \hspace{2cm}+
\int_\Omega\int_{Y\setminus D} 
\Cbb
\left(
\begin{array}{cc}
(Eu + E^y u^1)_{\alpha\beta} & (E^{yx}u)_{\alpha3}\\
(E^{yx}u)_{\alpha3} & 0
\end{array}
\right)
\\ \hspace{5cm}
\cdot
\left(
\begin{array}{cc}
(E\bar u+E^y\bar u^1)_{\alpha\beta}  &(E^{yx}\bar u)_{\alpha3}\\
(E^{yx}\bar u)_{\alpha3} & 0
\end{array}
\right)
\,dydx\\[10pt]
\displaystyle  
 \hspace{4cm}
=\int_\Omega \int_Yf_\alpha(x, y)\cdot \bar u_\alpha(x) 
+ f_3(x,y)\cdot \bar u_3(x,y) \,dydx,
\\[4ex]
\forall(\bar u_i,\bar u^1_\alpha,\bar \vartheta,\bar v_3,\bar w)\in  \mathcal{U}^h\times\mathcal{V}^h\times\mathcal{W}^h,
\end{array}
\right.
\end{equation}
where $\bar v_\alpha:=y_\alpha^R \bar \vartheta$.
\end{thm}

Remarks similar to those made in Remark \ref{onthemainresult} still hold.

The limit problem, even if it is not immediately apparent, is a second order problem since the term $(Eu)_{33}$ involves second derivatives of the function $u_\alpha$, as it can be easily deduced from the admissible set \eqref{defUh} . Thus we recover a second-order material as a limit of a first order material, see also \cite{BB2,PS}.

One can highlight the nonlocal effect in  the homogenized equation through the analogous of Theorem \ref{nonlocaleffect} as follows.

\begin{thm}\label{nonlocaleffecthomhom}
Let $(u_i,u^1_\alpha,\vartheta,v_3,w)$ be the solution of \eqref{Thelimitpbhom}. Let
$$E_f :=  \left(
\begin{array}{cc}
(E^y w)_{\alpha\beta} & (E^{yx}v)_{\alpha3}\\
(E^{yx}v)_{\alpha3} & (Eu)_{33}
\end{array}
\right),  \quad
 E_m := \left(
\begin{array}{cc}
(Eu + E^y u^1)_{\alpha\beta} & (E^{yx}u)_{\alpha3}\\
(E^{yx}u)_{\alpha3} & 0
\end{array}
\right).  $$
Then the sequence of solutions $u^\eps$ of \eqref{mainpbhomogenization} is such that
$$ Eu^\eps  \chi_{F_\eps}  \weak  \mint_{ D} E_f  \ dy \, \mbox{ in } \, L^2(\Omega;\Ro^{3 \times 3}), \quad  E  u^\eps \chi_{M_\eps}  \   \weak  \mint_{ Y \setminus D} E_m  \ dy \, \mbox{ in } \, L^2(\Omega;\Ro^{3 \times 3}), $$
$$ \eps u^\eps_\alpha \  \weak  u_\alpha \mbox{ in } \, H^1(\Omega),  $$
$$
 u^\eps_3 \  \weak  U(x) := \mint_Y u_3\ dy \mbox{ in } \, L^2(\Omega).
$$
$U$ may be written as 
$$U(x) = \mint_{D} u_3 \ dy + m_0(x) \mint_{Y \setminus D} f_3 \ dy + m_{00}(x),  $$ 
where $m_0$ and $m_{00}$ are given by
$$ m_0 (x) = \int_{Y \setminus D} z^0_3 \ dy, \quad   m_{00} (x) = \int_{Y \setminus D} z^{00}_3 \ dy,$$ being $z^0$ and $z^{00  }$ respectively the solutions of the following problems
\begin{equation}\label{elementarysolutionhom}
 \left\{ 
\begin{array}{ll}
z^0 \in L^\infty(\Omega;  \mathcal{Z}^{h}_0 ),   \\[2ex]   
 \displaystyle \int_{Y \setminus D} \Cbb(x,y)
\left(\begin{array}{cc}
(E^yz^0)_{\alpha\beta} & \frac 12 \partial^y_\alpha z^0_3\\
\frac 12 \partial^y_\alpha z^0_3 & 0
\end{array}\right)
\cdot
\left(\begin{array}{cc}
(E^y\bar z)_{\alpha\beta} & \frac 12 \partial^y_\alpha \bar z_3\\
\frac 12 \partial^y_\alpha \bar z_3 & 0
\end{array}\right)
\ dy
=   \\[4ex] \hspace{6cm}\displaystyle 
\int_{Y \setminus D} \bar z_3(y) \,dy\quad \mbox{for a.e.\ } x  \in \Omega,\\[2ex]     
\forall \bar z \in \mathcal{Z}^{h}_0 := \{ z : z_i \in H^1(Y), \, z_i=0\hbox{ in }   D \},
\end{array}
\right.
\end{equation}
\begin{equation}\label{secondsolutiononZhom}
 \left\{ 
\begin{array}{ll}
z^{00} \in L^2(\Omega; \mathcal{Z}^h_0),  \\[2ex]   
 \displaystyle \int_{\Omega \times(Y\setminus D)} \Cbb(x,y)
\left(\begin{array}{cc}
(E^yz^{00})_{\alpha\beta} & \frac 12 \partial^y_\alpha z^{00}_3\\
\frac 12 \partial^y_\alpha z^{00}_3 & 0
\end{array}\right)
\cdot
\left(\begin{array}{cc}
(E^y\bar z)_{\alpha\beta} & \frac 12 \partial^y_\alpha \bar z_3\\
\frac 12 \partial^y_\alpha \bar z_3 & 0
\end{array}\right)
 dx dy
\\[4ex] \hspace{2,5cm}\displaystyle =    
\int_{\Omega \times(Y\setminus D) } (f_3(x,y) - \mint_{Y\setminus D} f_3(x,y) \ dy )\bar z_3(x,y)) \,dx dy \\ 
\forall \bar z \in L^2(\Omega; \mathcal{Z}^h_0).
\end{array}
\right.
\end{equation}
\end{thm}
Similar remarks to those made in Remark \ref{onthenonlocaleffect} apply to the present contest after replacing the variable $x^\prime$ by $y$.

\section{\bf The rod problem: proofs of the main results}\label{rod_proofs}

We start by stating an a priori estimate and then work up to it.
\begin{thm}\label{comp_rod}
Let $u^\eps\in H^1_{dd}( \Omega)$ be the solution of problem \eqref{mainpb}.
Then
\begin{eqnarray}
&\displaystyle \sup_\eps \big(
\|E^\eps u^\eps\|_{L^2(F)}+\| \eps H^\eps u^\eps\|_{L^2(F)}+\|u^\eps\|_{H^1(F)}
\big)<+\infty,& \label{bduepsF}\\
& \displaystyle \sup_\eps \big(
\| \eps H^\eps u^\eps\|_{L^2(M)}+\|(u^\eps_1,u^\eps_2)\|_{H^1(M)}+\|u^\eps_3\|_{L^2(I;H^1(\omega\setminus \omin))}
\big)<+\infty,&\label{bduepsM}
\end{eqnarray}
in particular
\begin{equation}\label{bdueps}
\sup_\eps \big(
\| \eps H^\eps u^\eps\|_{L^2(\Omega)}+\|(u^\eps_1,u^\eps_2)\|_{H^1(\Omega)}+\|u^\eps_3\|_{L^2(I;H^1(\omega))}
\big)<+\infty.
\end{equation}
\end{thm}

To prove the theorem it will be convenient to set
$$
\mathcal{F}_\eps(v):= \int_\Omega (\chi_F+\eps^2\chi_M)\Cbb(x) E^\eps v\cdot E^\eps v\,dx.
$$
We now state and prove several short lemmas.

\begin{lem}\label{comp1}
For every $v\in H^1_{dd}(\Omega)$,  we have
$$
\mathcal{F}_\eps(v)\ge c\|E^\eps v\|^2_{L^2(F)}, \quad \mathcal{F}_\eps(v)\ge c\|v\|^2_{H^1(F)}.
$$
\end{lem}
{\bf { \proof}}
Using the strong ellipticity of  $\Cbb$, see \eqref{tensor} c),  we have
$$
\mathcal{F}_\eps(v)\ge m \|E^\eps v\|^2_{L^2(F)}\ge m \|E v\|^2_{L^2(F)}\ge c \|v\|^2_{H^1(F)},
$$
where the last inequality follows from Korn's inequality.
\QED

\begin{lem}\label{kornR}
Let $R,R^+,R^-\subset \Ro^3$ be three open regions with Lipschitz boundary such that $\bar R^+\cap \bar R^-$ has a strictly positive two dimensional measure and the closure of ${R^+\cup R^-}$
coincides with the closure of $R$, $\bar R$. Then, there exists a constant $c>0$ such that for every $v\in H^1(R;\Ro^3)$ we have that
$$
\|\nabla v\|_{L^2(R^+)}\le c\big( \|E v\|_{L^2(R^+)}+\|\nabla v\|_{L^2(R^-)} \big).
$$
\end{lem}
{\bf \proof}
Assume that such a constant $c$ does not exist. Then for every integer $n$ there exists a function
$v^n\in H^1(R;\Ro^3)$ such that
\begin{equation}\label{vn}
\|\nabla v^n\|_{L^2(R^+)}=1, \quad \|E v^n\|_{L^2(R^+)}+\|\nabla v^n\|_{L^2(R^-)} \le \frac 1n.
\end{equation}
Thus, $\|\nabla v^n\|_{L^2(R)}\le1$ and the sequence $\tilde v^n$ defined by

$$
\tilde v^n:=v^n-\mint_R v^n\,dx,
$$
 is bounded in $H^1(R;\Ro^3) $ so that up to a subsequence we have 
$$
\tilde v^n\weak \tilde v \quad \mbox{ in } H^1(R;\Ro^3),
$$
for some $\tilde v\in H^1(R;\Ro^3)$. 
But, from \eqref{vn} we infer that $\nabla \tilde v=0$ on $R^-$ and that
$E\tilde v=0$ on $R^+$.  Thus $\tilde v=k$  on $R^-$, with $k\in \Ro^3$.
As a consequence we first deduce that the trace of $\tilde v$ on  $\bar R^+\cap \bar R^-$
is $k$ and then that $\tilde v$ is equal to $k$ also on $R^+$, since $E\tilde v=0$ on $R^+$.
Thus $\tilde v=k$ on $R$. Since $\displaystyle \mint_R \tilde v^n\,dx=0$ we deduce that also
$\displaystyle \mint_R \tilde v\,dx=0$ and hence that $k=0,$ i.e., $\tilde v=0$ on $R$. 
By applying the standard Korn's inequality on $R^+$,
$$
\|\tilde v^n\|_{H^1(R^+)}\le c\big( \|E \tilde v^n\|_{L^2(R^+)}+\|\tilde v^n\|_{L^2(R^+)}\big)
$$
we see that $ \| \tilde v^n\|_{H^1(R^+)}$ approaches zero since the right-hand side converges to zero.
But this contradicts \eqref{vn}.
\QED

\begin{lem}\label{comp2}
There exists a constant $c>0$ such that for every $v\in H^1_{dd}(\Omega)$ we have
$$
\|H^\eps v\|_{L^2(M)} \le c \big( \|E^\eps v\|_{L^2(M)} + \|H^\eps v\|_{L^2(F)}\big).
$$
\end{lem}
{\bf \proof}
We first  extend $v$ by zero, without renaming it, to $\omega\times \Ro$.
Then we have that $ v\in H^1(\omega\times \Ro)$ and that, with the notation introduced in Remark \ref{omegaeps},
$\hat v\in H^1(\eps\omega\times \Ro)$.  We now partition $\eps\omega\times \Ro$
with slices of thickness $\eps$, that is we write
$$\eps\omega\times \Ro = \bigcup_{i \in I_\eps} S^i_\eps, \qquad S^i_\eps:=\eps\omega\times(\eps i, \eps (i+1)),$$
where $I_\eps$ is the set of the corresponding integer values of $i$. 
We also set
$$
F_\eps^i:= \eps \omin\times (\eps i, \eps (i+1)),\quad M_\eps^i:=S_\eps^i\setminus\bar F_\eps^i,
$$
and 
$$
\tilde S:=\omega\times (0,1),\quad \tilde F:=\omin\times (0,1),\quad\tilde M:=\tilde S\setminus \tilde F.
$$ 
Let $v^i:\tilde S\to \Ro^3$ be defined by
$$
v^i(x):=\hat v(\eps x_1,\eps x_2, \eps (i+x_3)).
$$
By Lemma \ref{kornR}, there exists a positive constant $c$ such that  for every $i$,
$$
\int_{\tilde M}|\nabla v^i|^2\,dx \le c 
\Big(
\int_{\tilde M}|E v^i|^2\,dx+\int_{\tilde F}|\nabla  v^i|^2\,dx
\Big),
$$
and changing variables, we find
$$
\int_{M^i_\eps}|\nabla \hat v|^2\,dx \le c 
\Big(
\int_{M^i_\eps}|E \hat v|^2\,dx+\int_{F^i_\eps}|\nabla \hat v|^2\,dx
\Big),
$$
where the constant $c$ does not depend neither on $i$ nor on $\eps$. Summing over $i$ we find
$$
\int_{M_\eps}|\nabla \hat v|^2\,dx \le c 
\Big(
\int_{M_\eps}|E \hat v|^2\,dx+\int_{F_\eps}|\nabla \hat v|^2\,dx
\Big),
$$
and changing variables, using again the notation introduced in Remark \ref{omegaeps}, we conclude the proof.
\QED

\begin{lem}\label{comp3}
There exists a constant $c>0$ such that for every $v\in H^1_{dd}(F)$ we have
$$
\|H^\eps v\|_{L^2(F)}+ \|(R^\eps)^{-1}v \|_{L^2(F)} \le \frac c\eps  \|E^\eps v\|_{L^2(F)}.
$$
\end{lem}
{\bf \proof}
This follows immediately from a rescaled Korn's inequality proved in \cite{ABP94}. Indeed, using the notation introduced in Remark \ref{omegaeps}, let $\hat v=(R^\eps)^{-1} v\circ ( r^\eps)^{-1}$. Then in \cite{ABP94} it has been proved that
$$
\int_{F_\eps} |\nabla\hat v|^2+|\hat v|^2\,dx \le \frac c{\eps^2} \int_{F_\eps} |E\hat v|^2\,dx.
$$
The result follows by simply changing variables.
\QED

\begin{lem}\label{comp4}
For every $v\in H^1_{dd}(\Omega)$ we have
$$
\mathcal{F}_\eps(v)\ge c\| \eps H^\eps v\|^2_{L^2(M)}.
$$
\end{lem}
{\bf \proof}
Since $\Cbb$ is uniformly positive definite we have
\begin{eqnarray*}
\mathcal{F}_\eps(v)&\ge& m\big(\|E^\eps v\|^2_{L^2(F)}+\eps^2 \|E^\eps v\|^2_{L^2(M)}\big)
\\
&\ge& c\eps^2\big(\|H^\eps v\|^2_{L^2(F)}+ \|E^\eps v\|^2_{L^2(M)}\big),\\
&\ge&c\eps^2\|H^\eps v\|^2_{L^2(M)}
\end{eqnarray*}
where we first used Lemma \ref{comp3} and then Lemma \ref{comp2}.
\QED

\begin{lem}\label{comp5}
For every $v\in H^1_{dd}(\Omega)$ we have
$$
\mathcal{F}_\eps(v)\ge c\big(\|v'\|^2_{H^1(M)}+\|v_3\|^2_{L^2(I;H^1(\omega\setminus \omin))}\big),
$$
with $v'=(v_1,v_2)$.
\end{lem}
{\bf \proof}
From Lemma \ref{comp4},  \eqref{components},  Poincar\'e inequality, and small $\eps$,  we have  that
$$
\mathcal{F}_\eps(v)\ge c\big(\| \frac{\nabla' v'}{\eps}\|^2_{L^2(M)}+\|{\partial_3 v'}\|^2_{L^2(M)}\big)\ge c \|\nabla v'\|^2_{L^2(M)}\ge c \| v'\|^2_{H^1(M)}.
$$
Also, from Lemma \ref{comp4} and \eqref{components} we find that
\begin{equation}\label{bddv3}
\mathcal{F}_\eps(v)\ge c\|{\nabla' v_3}\|^2_{L^2(M)}.
\end{equation}
From the following Poincare's type inequality
$$
\int_M |v_3|^2\,dx\le c\big( \int_\Omega |\nabla'v_3|^2\,dx+\int_F |v_3|^2\,dx\big),
$$
Lemma \ref{comp1} and \eqref{bddv3}, it follows that
$$
\mathcal{F}_\eps(v)\ge c\|{ v_3}\|^2_{L^2(M)},
$$
and hence the statement of the Lemma follows.
\QED

\begin{lem}\label{comp6}
Let $u^\eps\in H^1_{dd}( \Omega)$ be the solution of problem \eqref{mainpb}.
Then
$$
\sup_\eps \mathcal{F}_\eps(u^\eps)<+\infty.
$$
\end{lem}
{\bf \proof}
By taking $\varphi=u^\eps$ in problem \eqref{mainpb}, we find
$$
\mathcal{F}_\eps(u^\eps) \le \|f\|_{L^2(\Omega)}\|u^\eps\|_{L^2(\Omega)},
$$
while from Lemma \ref{comp1} and Lemma \ref{comp5} we deduce that
$$
c \|u^\eps\|^2_{L^2(\Omega)}\le \mathcal{F}_\eps(u^\eps).
$$
The Lemma follows from the previous two inequalities.
\QED

From Lemmas \ref{comp1}, \ref{comp3},  \ref{comp4}, \ref{comp5}  and \ref{comp6}, one deduces
\eqref{bduepsF}, \eqref{bduepsM}, and \eqref{bdueps}; thus, the proof of Theorem \ref{comp_rod} is achieved.

We now prove some compactness properties using the previous a priori estimates.
 
\begin{lem}\label{convu3}
There exist $(u, v, w) \in \mathcal{U} \times \mathcal{V} \times \mathcal{W}$ and $z_3\in L^2(I;H^1(\omega))$ with
$z_3=0$ almost everywhere in $F$, such that the sequence of solutions $u^\eps$ of problem \eqref{mainpb} fullfills (up to a subsequence)
$$
u^\eps_\alpha\weak u_\alpha \mbox{ in }H^1(\Omega),
$$
$$
u^\eps_3\weak  u_3+z_3 \mbox{ in }L^2(I;H^1(\omega)).
$$
Moreover,
$$
(E^\eps u^\eps)_{33}\chi_F\weak (Eu)_{33}\chi_F, \quad
\eps(E^\eps u^\eps)_{33}\chi_M\weak 0,\quad \eps(E^\eps u^\eps)_{\alpha3}\chi_M\weak \frac 12 \partial_\alpha z_3\chi_M,
$$ 
$$
(E^\eps u^\eps)_{\alpha3}\chi_F\weak (E v)_{\alpha3}\chi_F, \quad
(E^\eps u^\eps)_{\alpha\beta}\chi_F\weak (E w)_{\alpha\beta}\chi_F  \mbox{ in } L^2(\Omega),
$$
in $L^2(\Omega)$. 
\end{lem}

{\bf \proof}
From \eqref{bdueps}, up to subsequences, it follows that
$$
u^\eps_\alpha\weak u_\alpha \mbox{ in }H^1(\Omega),
\qquad
u^\eps_3\weak \tilde u_3 \mbox{ in }L^2(I;H^1(\omega)),
$$
for some $ u_\alpha \in H^1_{dd}(\Omega)$ and $\tilde u_3\in L^2(I;H^1(\omega))$.
From \eqref{components} and \eqref{bdueps} we deduce that $\nabla' u^\eps_\alpha\to 0$ in 
$L^2(\Omega;\Ro^2)$ thus $\nabla' u_\alpha=0$. Hence there exist $\xi_\alpha\in H^1_{0}(0,\ell)$
such that
$$
u_\alpha(x)=\xi_\alpha(x_3) \mbox{ for a.e. }x\in \Omega .
$$
From \eqref{bduepsF} it follows that $u^\eps_3\weak \tilde u_3 \mbox{ in }H^1(F)$ and since
$(Eu^\eps)_{i\alpha}\to 0$ in $L^2(F)$, we have that $(u_1,u_2,\tilde u_3)\in \mathcal{BN}_{dd}(F)$, that is: $\xi_\alpha\in H^2_{0}(0,\ell)$ and there exists a function $\xi_3\in H^1_{0}(0,\ell)$ such that
$$
u_\alpha=\xi_\alpha, \quad \tilde u_3=\xi_3-x_\alpha \xi_\alpha'\quad \mbox{ in }F.
$$
Set 
$$
u_3:=\xi_3-x_\alpha \xi_\alpha'\in H_{dd}^1(\Omega), \quad z_3:=\tilde u_3-u_3\in L^2(I;H^1(\omega)).
$$
Then, $z_3=0$ in $F$, $u:=(u_1,u_2,u_3)\in \mathcal{BN}_{dd}(\Omega)$,
and
$
u^\eps_3\weak u_3+z_3 \mbox{ in }L^2(I;H^1(\omega)).
$

Since $(E^\eps u^\eps)_{33}=\partial_3 u^\eps_3$ we have that
$$
(E^\eps u^\eps)_{33}\chi_F\weak \partial_3\tilde u_3\chi_F= \partial_3 u_3 \chi_F=(Eu)_{33}\chi_F
$$
in $L^2(\Omega)$. Since $\eps(E^\eps u^\eps)_{33}=\eps(H^\eps u^\eps)_{33}$ from \eqref{bduepsM} it follows that 
$\eps(E^\eps u^\eps)_{33}$ is weakly convergent in $L^2(M)$, but since $\eps u^\eps_3\to 0$ in $L^2(M)$ we deduce that 
$$
\eps(E^\eps u^\eps)_{33}\chi_M\weak 0 \mbox{ in }L^2(\Omega).
$$
Finally,
$$
\eps(E^\eps u^\eps)_{\alpha3}\chi_M=\frac 12 (\partial_\alpha u^\eps_3+\partial_3 u^\eps_\alpha)\chi_M\weak  \frac 12 (\partial_\alpha u_3+\partial_\alpha z_3+\partial_3 u_\alpha)\chi_M=
(Eu)_{\alpha3}+\frac 12 \partial_\alpha z_3
$$
in $L^2(\Omega)$. The stated convergence result follows by recalling that $u\in \mathcal{BN}_{dd}(\Omega)$ and hence that $(Eu)_{\alpha3}=0$.

For the proof of the last two convergences, see  \cite{MS1} and \cite{CS}.
\QED

In  Lemmas 3.1 and 3.2 of \cite{FMP2010}  it is shown that
\begin{equation}\label{H21}
\eps(H^\eps u^\eps)_{21}\chi_F\weak \vartheta \chi_F\mbox{ in }L^2(\Omega), 
\end{equation}
where $ \vartheta \in H_0^1(0,\ell)$; moreover, with $v$ as in Lemma \ref{convu3} we have $ (Ev)_{\alpha3}= \frac12( \vartheta' x_\alpha^R + \partial_\alpha v_3)$.
 
\begin{lem}\label{convza}
Let $\vartheta \in H^1_0(I)$ be as in \eqref{H21}. There exists $z_\alpha\in L^2(I;H^1(\omega))$ with
$z_\alpha=0$ almost everywhere in $F$, such that up to a subsequence,
$$
\frac{u^\eps_\alpha}{\eps}-\mint_\omin \frac{u^\eps_\alpha}{\eps}\,dx'\weak  \vartheta x^R_\alpha+z_\alpha \mbox{ in }L^2(I;H^1(\omega)).
$$
Moreover,
$$
\eps(E^\eps u^\eps)_{\alpha\beta}\chi_M\weak (Ez)_{\alpha\beta}\chi_M
\mbox{ in } L^2(\Omega).$$
\end{lem}

{\bf \proof}
Let
$$
z^\eps_\alpha:= \displaystyle \frac{u^\eps_\alpha}{\eps}-\mint_\omin \frac{u^\eps_\alpha}{\eps}\,dx'
$$
Since $ \displaystyle \mint_\omin {z^\eps_\alpha} \, dx'=0$,  by Poincare's inequality we have 

$$
\int_\omega |z^\eps_\alpha|^2\,dx'\le C \int_\omega |\nabla' z^\eps_\alpha|^2\,dx'
=C \frac 1{\eps^2}\int_\omega |\nabla' u^\eps_\alpha|^2\,dx'\le C \int_\omega |\eps H^\eps u^\eps|^2\,dx',
$$
a.e.\ in $I =  (0 ,\ell)$. Integrating in $x_3$ over $I$ we deduce, thanks to \eqref{bdueps},  that $z^\eps_\alpha$ is bounded in $L^2(I;H^1(\omega))$  so that up to a subsequence,
$$
z^\eps_\alpha\weak \tilde z_\alpha \mbox{ in }L^2(I;H^1(\omega)),
$$
for some $\tilde z_\alpha \in L^2(I;H^1(\omega))$ such that $\displaystyle \mint_\omin {\tilde z_\alpha}\,dx'=0$.
From \eqref{bduepsF} we have that $\eps E^\eps u^\eps\to 0$ in $L^2(F;\Ro^{3\times3})$ and hence, taking into account \eqref{H21}, we have that
$$
\eps (H^\eps u^\eps)_{\alpha\beta}\weak 
\left( 
\begin{array}{cc}
0&-\vartheta\\
\vartheta&0
\end{array}
\right)  \mbox{ in } L^2(F;\Ro^{2\times2}).
$$
But $\eps (H^\eps u^\eps)_{\alpha\beta}=\frac 1\eps \partial_\beta u^\eps_\alpha= \partial_\beta z^\eps_\alpha$ and hence
$$
\partial_1 \tilde z_1=\partial_2 \tilde z_2=0, \qquad \partial_1 \tilde z_2=-\partial_2 \tilde z_1=\vartheta\qquad \mbox{ in } F.
$$
By integration we find
$$
\tilde z_\alpha= \vartheta x_\alpha^R +d_\alpha, \quad \mbox{inÊ} \ F,
$$
with $d_\alpha\in L^2(I)$. But since $\displaystyle \mint_\omin {\tilde z_\alpha}\,dx'=0$,  we deduce, by invoking \eqref{centermass}, that $d_\alpha=0$.

We set
$$
z_\alpha:=\tilde z_\alpha- \vartheta x_\alpha^R \, \, \mbox{in} \,  \Omega.
$$
Then $z_\alpha=0$ in $F$  and $z_\alpha\in L^2(I;H^1(\omega))$ . Also,
$$
\eps (E^\eps u^\eps)_{\alpha\beta}\chi_M=\frac{1}{2\eps}(\partial_\alpha u^\eps_\beta+
\partial_\beta u^\eps_\alpha)\chi_M=(Ez^\eps)_{\alpha\beta}\chi_M\weak
(E\tilde z)_{\alpha\beta}\chi_M=(E z)_{\alpha\beta}\chi_M
$$
in $L^2(\Omega)$.
\QED

 \noindent
 {\bf Proof of Theorem \ref{Mainresult}.} 
 The convergences stated in Theorem \ref{Mainresult} with strong convergence replaced by weak convergence are proved in Lemmas \ref{convu3} and \ref{convza}. At the end of the proof we will show that indeed are strong. To identify the limit problem take test function $\varphi$ in \eqref{mainpb} of the following form 
\begin{equation}\label{testPhi}
\left\{\begin{array}{ll}
& \varphi_\alpha(x)=\bar u_\alpha(x_3)+ \eps \bar v_\alpha(x)+ \eps \bar z_\alpha(x)+ \eps^2 \bar w_\alpha(x)\\
& \varphi_3(x)=\bar u_3(x)+  \bar z_3(x)+ \eps \bar v_3(x),
\end{array}\right.
\end{equation}
with
$$
\bar u=(\bar u_\alpha,\bar u_3)\in\mathcal{BN}_{dd}(\Omega),\quad
\bar v=(\bar v_\alpha,\bar v_3)\in\mathcal{R}_{dd}(\Omega)\cap H^1_{dd}(\Omega)
$$
and
$$
\bar w=(\bar w_\alpha,0)\in \mathcal{RD}^\perp_{2}(F)\cap H^1_{dd}(\Omega),\quad
\bar z=(\bar z_\alpha,\bar z_3)\in H^1_{dd}(\Omega), \mbox{ with }\bar z_i=0 \mbox{ in }F.
$$
Then
\begin{eqnarray*}
(E^\eps\varphi)_{\alpha\beta}&=&\frac 1\eps (E\bar z)_{\alpha\beta}\chi_M+(E\bar w)_{\alpha\beta},\\
(E^\eps\varphi)_{\alpha3}&=&\frac 1{2\eps} \partial_\alpha\bar z_3\chi_M+
\frac 1{2} \partial_3\bar z_\alpha\chi_M+(E\bar v)_{\alpha3}+\eps
(E\bar w)_{\alpha 3},\\
(E^\eps\varphi)_{33}&=&(E\bar u)_{33}+(E\bar z)_{33}\chi_M + \eps
(E\bar v)_{\alpha 3}.
\end{eqnarray*}
Taking $\varphi$ in \eqref{mainpb} as in \eqref{testPhi}, letting $\eps$ go to zero, and using the density of $\mathcal{R}_{dd}(\Omega)\cap H^1_{dd}(\Omega)$,  $\mathcal{RD}^\perp_{2}(F)\cap H^1_{dd}(\Omega)$ and $\mathcal{Z} \cap H^1_{dd}(\Omega)$  into $\mathcal{V}$, $\mathcal{W}$ and $\mathcal{Z}$ respectively, we deduce \eqref{Thelimitpb}.

We now prove the strong convergences stated in the theorem. Set
$$E_f :=  \left(\begin{array}{cc}
(Ew)_{\alpha\beta} & (Ev)_{\alpha 3}\\
(Ev)_{\alpha 3} & (Eu)_{3 3}
\end{array}\right),  \quad
 E_m := \left(\begin{array}{cc}
(Ez)_{\alpha\beta} & \frac 12 \partial_\alpha z_3\\
\frac 12 \partial_\alpha z_3 & 0
\end{array}\right).  $$
Taking $\varphi=u^\eps$ in \eqref{mainpb} and passing to the limit we find
\begin{eqnarray*}
\lim_{\eps\to 0} \int_\Omega (\chi_F+\eps^2\chi_M)\Cbb E^\eps u^\eps\cdot E^\eps u^\eps\,dx&=&\lim_{\eps\to 0}
\int_\Omega f \cdot u^\eps \,dx\\
&=& \int_\Omega f_\alpha\bar u_\alpha+f_3(\bar u_3+\bar z_3)\,dx\\
&=&\int_\Omega \Cbb E_f\cdot E_f\chi_F+\Cbb E_m\cdot E_m\chi_M\,dx.
\end{eqnarray*}
Since the norm $A\mapsto \int_\Omega \Cbb A\cdot A\, dx$ is equivalent to $A\mapsto \|A\|_{L^2(\Omega)}$, it follows that
$$
 \lim_{\eps\to 0}\|(\chi_F+\eps\chi_M) E^\eps u^\eps\|_{L^2(\Omega)}=
 \|\chi_F E_f+\chi_M E_m\|_{L^2(\Omega)},
$$
which, together with the weak convergence, implies that
\begin{equation}\label{str}
(\chi_F+\eps\chi_M) E^\eps u^\eps\to \chi_F E_f+\chi_M E_m\quad\mbox{in }L^2(\Omega).
\end{equation}
We now prove the strong convergence of the displacement components.
By Korn inequality  we have, for $\eps$ small,  that $\|E^\eps u^\eps\|_{L^2(F)}\ge  \|E u^\eps\|_{L^2(F)}\ge c \|u^\eps\|_{H^1(F)}$ and hence, by means of \eqref{str}, we deduce that $u^\eps$ is a Cauchy sequence in $H^1(F)$. By means of Lemma \eqref{comp3} first and then Lemma \eqref{comp2} we deduce that 
$\eps H^\eps u^\eps$ is a Cauchy sequence in $L^2(\Omega)$ and then, with the same argument used in the proof of Lemma \eqref{comp5}, we  find that $u^\eps_\alpha$ is a Cauchy sequence in $H^1(M)$ and that $u^\eps_3$ is a Cauchy sequence in $L^2(I;H^1(\omega\setminus\omin))$. This completes the proof since 
we have already shown that $u^\eps$ is a Cauchy sequence in $H^1(F)$. 
\QED

%%%%%%%%%%%%%%%%

\noindent
{\bf Proof of Theorem \ref{nonlocaleffect}.}
All the convergences stated in Theorem \ref{nonlocaleffect}  may be obtained choosing a test function $\varphi \in  L^2(\Omega)$ that depends only on the $x_3$ variable in the  corresponding convergences stated in  Theorem \ref{Mainresult}. To prove that $U$ can be written as
$$
U(x_3) = \mint_\omega u_3 \ dx' + m_0(x_3) \mint_{\omega \setminus \omin} f_3 \ dx' + m_{00}(x_3),  
$$  
we first recall that, by definition,  $\displaystyle U(x_3)=\mint_\omega u_3 +z_3\ dx'$, where $z_3$ is the solution of problem \eqref{Thelimitpb} with $ \bar u = \bar v = \bar w = 0$, that is
\begin{equation}\label{limpbonZ}
\left\{
\begin{array}{ll}
\displaystyle z_3 \in  \mathcal{Z}, \\[1ex]
\displaystyle    
\int_\Omega \Cbb(x)
\left(\begin{array}{cc}
(Ez_3)_{\alpha\beta} & \frac 12 \partial_\alpha z_3\\
\frac 12 \partial_\alpha z_3 & 0
\end{array}\right)
\cdot
\left(\begin{array}{cc}
(E\bar z)_{\alpha\beta} & \frac 12 \partial_\alpha \bar z_3\\
\frac 12 \partial_\alpha \bar z_3 & 0
\end{array}\right)
\chi_M\,dx
= \\[4ex]\hspace{8cm}
\displaystyle  
\int_\Omega f_3(x)\bar z_3(x) \,dx,\\
  \forall \, \bar z_3 \in  \mathcal{Z}.
\end{array}
\right.
\end{equation}
Since problem \eqref{limpbonZ} is linear we can write $z_3$ as $z_3 = z^1 + z^{00}$, where
$z^{00}$ is the solution of problem \eqref{limpbonZ} but with $f_3(x)$ replaced by 
$f_3(x)- \displaystyle \mint_{\omega \setminus \omin} f_3(x',x_3) \ dx'$, and $z^1$ is the solution
of problem \eqref{limpbonZ} but with $f_3(x)$ replaced by 
$\displaystyle \mint_{\omega \setminus \omin} f_3(x',x_3) \ dx'$. By linearity and uniqueness one can check that 
$z^1= \displaystyle \mint_{\omega \setminus \omin} f_3(x',x_3) \ dx' z^0$, where $z^0$ is the solution of problem \eqref{elementarysolution}. The representation of $U$ then immediately follows. 
\QED

\section{\bf The homogenization problem: proofs of the main results}\label{homo_proofs}

We start by proving some a priori estimates.

\begin{lem}\label{comp7}
The sequence of solutions $u^\eps$ of problem \eqref{mainpbhomogenization} satisfies the following a priori estimates
\begin{eqnarray}
&  \displaystyle \sup_\eps \|  u_3^\eps\|_{L^2(\Omega)}<+\infty,  & \label{bdh1} \\    
& \displaystyle \sup_\eps \| \eps u^\eps\|_{H_{dd}^1(\Omega)}<+\infty. & \label{bdh2} 
\end{eqnarray}
Moreover, setting 
$$
\mathcal{J}_\eps(u^\eps):= \int_\Omega ( \chi_{F_\eps} + \eps^2  \chi_{M_\eps}) \vert E u^\eps \vert^2 \ dx,
$$
we also have
\begin{equation}
\sup_\eps \mathcal{J}_\eps(u^\eps) <+\infty.  \label{bdh3} 
\end{equation}
\end{lem}

{\bf \proof}
By Korn inequality, there exists a constant $C$ such that
$$\int_I \int_D \vert v \vert^2 \ dy dx_3 \leq C \int_I \int_D \vert E v \vert^2 \ dy dx_3 \quad \forall \, v \in H^1_{dd}( D \times I).     $$
With $v$ as  
$$v_\alpha(y,x_3) = \eps u^\eps_\alpha (\eps y + \eps i, x_3); \quad v_3(y,x_3) =  u^\eps_3 (\eps y + \eps i, x_3), \quad \forall \, y \in D \times I, $$
we obtain, after the change of variable $x:= (\eps y + \eps i, x_3) \in D_\eps^i \times I := F_\eps^i$, that
$$
\begin{array}{ll}
 \displaystyle  \int_{ÊF_\eps^i} ( \eps^2 \vert u_\alpha^\eps \vert^2   + \vert u_3^\eps \vert^2 ) \ dx & \leq C \displaystyle \int_{ÊF_\eps^i} \sum_{\alpha \beta} ( \eps^4\vert (E u^\eps)_{\alpha \beta} \vert^2  +  \eps^2\vert (E u^\eps)_{\alpha 3} \vert^2)  \ + \vert (E u^\eps)_{33} \vert^2    dx  \\
& \displaystyle   \leq C  \int_{ÊF_\eps^i}  \vert E u^\eps \vert^2  \ dx.
\end{array} 
$$
Taking the sum over $ i$ in the previous inequality, we get the estimate
\begin{eqnarray}\label{bdF}
\displaystyle \int_{F_\eps} \vert u^\eps_3 \vert ^2  \ dx \leq C \int_{F_\eps} \vert Eu^\eps \vert ^2  \ dx.
\end{eqnarray}
We now apply the following Poincar\'e 's type inequality 
$$
\int_M |v|^2\,dx\le c\big( \int_\Omega |\nabla'v |^2\,dx+\int_F |v|^2\,dx\big),
$$
with $ \Omega : = Y\times I, \, \, F:= D \times I, \, \, M:= (Y \setminus \bar D) \times I$ and
 $ v(y,x_3) = u^\eps_3 (\eps y + \eps i, x_3)$ so that using the same change of variable as above and the classical Korn inequality $\|  \eps u^\eps \|_{H^1(\Omega)}  \leq C \|  \eps E u^\eps \|_{L^2(\Omega)} $, we infer, after taking the sum over $i$, that
\begin{eqnarray}\label{bdM}
\displaystyle \int_{M_\eps} \vert u^\eps_3 \vert ^2  \ dx \leq C \int_{\Omega} \vert  \eps E u^\eps \vert ^2  + \int_{F_\eps} \vert u^\eps_3 \vert ^2   \ dx.
\end{eqnarray}
Since by definition $ \mathcal{J}_\eps(u^\eps)=  \displaystyle \int_\Omega ( \chi_{F_\eps} + \eps^2  \chi_{M_\eps}) \vert E u^\eps \vert^2 \ dx $, 
we deduce from \eqref{bdF}  and \eqref{bdM} that
\begin{eqnarray}\label{bd3}
\displaystyle \int_\Omega \vert u^\eps_3 \vert ^2  \ dx \leq C\mathcal{J}_\eps(u^\eps). 
\end{eqnarray}
We now take $\varphi = u^\eps$ in \eqref{mainpbhomogenization}. Applying the Cauchy-Schwarz inequality in the right hand side and using the strong ellipticity of $ \mathbb{C}$ together with \eqref{bd3} and the Korn inequality $\|  \eps u^\eps \|_{H^1(\Omega)}  \leq C \|  \eps E u^\eps \|_{L^2(\Omega)} $, we get
\begin{equation}\label{finalbound}
\displaystyle  \mathcal{J}_\eps(u^\eps) 
\leq c \parallel \eps u^\eps_\alpha \parallel_{L^2(\Omega)} + c  \parallel  u^\eps_3 \parallel_{L^2(\Omega)} 
\leq c \parallel \eps u^\eps_\alpha \parallel_{H^1(\Omega)} + c  \parallel  u^\eps_3 \parallel_{L^2(\Omega)}
\leq c \sqrt{  \mathcal{J}_\eps(u^\eps)  }.
\end{equation}
From \eqref{finalbound}, we first deduce  \eqref{bdh3} and then, by means of \eqref{finalbound}, we deduce \eqref{bdh1} and \eqref{bdh2}.
\QED
\begin{lem}\label{lem412}
There exist
$(u_i,u^1_\alpha,\vartheta) \in \mathcal{U}^h$,  $v_3 \in \mathcal{V}^h$, and $w \in \mathcal{W}^h$ such that the sequence of solutions $u^\eps$ of problem \eqref{mainpbhomogenization} fullfills (up to a subsequence)
 \begin{eqnarray}
&  \displaystyle   \eps u^\eps_\alpha \weak u_\alpha \quad \mbox{in }  H^1_{dd}(\Omega),  & \label{tsvh100} \\  
& \displaystyle   \eps \nabla' u^\eps_\alpha \twoscale \nabla' u_\alpha +  \nabla_y' u^1_\alpha, & \label{tsvh20} \\
& u^\eps_3 \twoscale  u_3, & \label{tsvh40} \\
&\eps \nabla' u^\eps_3 \twoscale  \nabla'_y  u_3.& \label{tsvh50} 
 \end{eqnarray}
Moreover, setting $v_\alpha:=y_\alpha^R \vartheta$, we have
$$
(E u^\eps)\chi_{F_\eps} \twoscale 
\left(
\begin{array}{cc}
(E^y w)_{\alpha\beta} & (E^{yx}v)_{\alpha3}\\
(E^{yx}v)_{\alpha3} & (Eu)_{33}
\end{array}
\right)
\chi_D(y),
$$ 
$$ \eps(E u^\eps)\chi_{M_\eps} \twoscale  
\left(
\begin{array}{cc}
(Eu + E^y u^1)_{\alpha\beta} & (E^{yx}u)_{\alpha3}\\
(E^{yx}u)_{\alpha3} & 0
\end{array}
\right)
\chi_{Y \setminus D} (y),  
$$ 
where the operator $E^{yx}$ has been defined in \eqref{Eyx}.
\end{lem}

{\bf \proof}
By \eqref{bdh1}, \eqref{bdh2}, and  classical compactness results related to two-scale convergence, recalled in Section \ref{problem}, there exist   $ u_\alpha \in H^1_{dd} (\Omega)$, $ u^1_\alpha, u_3 \in L^2(\Omega; H^1_{\#}(Y)) $ such that, up to subsequences, 
the convergences stated in \eqref{tsvh100}--\eqref{tsvh50} hold.
 From \eqref{bdh3} it follows that 
 \begin{equation}\label{Ebd}
 \sup_\eps \|\chi_{F_\eps}Eu^\eps \|_{L^2(\Omega)}<+\infty\qquad
 \sup_\eps \|{\eps}Eu^\eps \|_{L^2(\Omega)}<+\infty.
 \end{equation}
 Thus, up to subsequences, we have that
 $$\chi_{F_\eps}\eps(Eu^\eps)_{\alpha\beta}\twoscale 0,$$
 while taking into account \eqref{tsvh20} it also follows that
 $$
 \chi_{F_\eps}\eps(Eu^\eps)_{\alpha\beta}\twoscale ((Eu)_{\alpha\beta} + (E^yu^1)_{\alpha\beta})\chi_D(y). 
 $$
 Hence
 $$
 (Eu)_{\alpha\beta} + (E^yu^1)_{\alpha\beta}=0 \mbox{ for a.e.}\
(x,y)\in\Omega\times D.
$$
Since $u$ does not depend on $y$, from the above equation we deduce that
$$
u^1_\alpha(x,y)=-(Eu(x))_{\alpha\beta}y_\beta+\tilde \vartheta(x) y^R_\alpha\in L^2(\Omega;H^1_{m}(D))
$$
for some function $\tilde \vartheta \in L^2(\Omega)$.
We remark that for every $x$ the sum $\tilde \vartheta(x) y^R_\alpha$ denotes an infinitesimal rigid rotation. Since
$$
(Eu)_{\alpha\beta}y_\beta=\partial_\beta u_\alpha y_\beta + \frac{\partial_2 u_1-\partial_1 u_2}2 y^R_\alpha,
$$
we can write
$$
u^1_\alpha(x,y)=- \partial_\beta u_\alpha (x) y_\beta+ \vartheta(x) y^R_\alpha
\mbox{ for a.e.}\ (x,y)\in\Omega\times D,
$$
where we set
$$
\vartheta:=\tilde \vartheta-\frac{\partial_2 u_1-\partial_1 u_2}2.
$$

By considering the component $(Eu^\eps)_{33}$, and \eqref{Ebd} we get that 
$$
\dis \partial_3 u^\eps_3  \chi_{F_\eps} \twoscale K_{33}(x,y) \chi_D(y),
$$ 
for some function $K_{33}\in L^2(\Omega\times D)$.
For $\phi \in \mathcal{D}(\Omega)$ and $\psi \in C_{\#}(Y)$, we have
\begin{equation}\label{cvu300}
\dis \int_\Omega \partial_3 u^\eps_3 \chi_{F_\eps}  \phi(x) \psi(\frac{x'}{\eps}) \ dx = - \int_\Omega u^\eps_3 \chi_{F_\eps}  \partial_3 \phi(x) \psi(\frac{x'}{\eps}) \ dx 
\end{equation}
and by passing to the limit we find
$$
\int_\Omega \int_D K_{33} (x,y)  \phi(x) \psi(y) \ dx dy
=
-\int_\Omega \int_D  u_3 (x,y) \partial_3 \phi(x) \psi(y) \ dx dy.
$$
Therefore, $u_3  \in L^2(\omega \times D; H^1(I))$ and 
 $ \dis \partial_3 u_3 =  K_{33} \, \, \mbox{a.e.\ in } \Omega \times D.$
In addition, we have $  u_3  \in L^2(\omega \times D; H_0^1(I))$. Indeed, identity \eqref{cvu300} holds also for test functions  $\phi \in C^\infty(\bar \Omega)$ and  $\psi \in C_{\#}(Y)$, since $u^\eps_3\in H^1_{dd}(\Omega)$, and hence we  deduce that
\begin{eqnarray*}  
\int_\Omega \int_Y \partial_3  u_3  \chi_D(y)\phi(x) \psi(y) \ dx dy
%\label{cvu4}
=
-\int_\Omega \int_Y  u_3 (x,y) \chi_D(y) \partial_3 \phi(x) \psi(y) \ dx dy.
\end{eqnarray*}
Integrating by parts the integral on the right, we conclude that 
\begin{equation}\label{u30l}
u_3(x', 0,y)  =   u_3(x', \ell,y) = 0.
\end{equation}
From  \eqref{Ebd} we deduce that
$$
2 \eps (E u^\eps)_{\alpha3} \chi_{F_\eps} \twoscale 0, 
$$
while, by means of \eqref{tsvh100} and \eqref{tsvh50} we conclude that
$$
2 \eps (E u^\eps)_{\alpha3} \chi_{F_\eps} =  ( \eps \partial_\alpha u^\eps_3 +  \eps \partial_3 u^\eps_\alpha) \chi_{F_\eps}\twoscale (\partial^y_\alpha u_3  + \partial_3u_\alpha) \chi_D(y).
$$
 In the above limit we used the fact that $\eps\partial_3 u^\eps_\alpha\twoscale \partial_3 u_\alpha$ since there is no fast variable in the $x_3$-direction. Thus
 $$ 
 (\partial^y_\alpha u_3  + \partial_3u_\alpha) \chi_D(y) = 0,
 $$
 from which we easily deduce that 
$$ \dis u_3 (x,y) = - y_\alpha \partial_3u_\alpha(x) + \hat u_3 (x), \, \mbox{ for a.e. } (x,y)\in \Omega \times D,$$
for some $\hat u_3\in  L^2(\omega;H^1_0(I)).$
Note that this equality implies that $u_\alpha$ also belongs to the space $L^2(\omega;H_0^2(I))$.

 We now prove the convergent of the rescaled strains.
 Thanks to \eqref{bdh3}, there exist symmetric matrix fields $K, G \in L^2(\Omega;\Ro^{3 \times 3})$ such that, up to subsequences, 
 \begin{eqnarray}
& \displaystyle E u^\eps \chi_{F_\eps}  \twoscale K \chi_D(y), & \label{tsvh6} \\
&  \displaystyle \eps E u^\eps \chi_{M_\eps}  \twoscale G \chi_{Y \setminus D}(y). & \label{tsvh7}
 \end{eqnarray}
 
 We now characterize the components of $K$ and $Z$.
That $K_{33}=\partial_3  u_3 $ a.e.\ in $\Omega\times D$ has been already proven above.
The proof of the characterization of the two-scale limits $K_{\alpha3} =\frac 12 (y^R_\alpha \partial_3 \hat\vartheta+\partial^y_\alpha v_3)$ and $ K_{\alpha\beta} =  (E^y w)_{\alpha\beta}$ for  some $v_3 \in \mathcal{V}^h$, $\hat\vartheta\in L^2(\omega;H^1_0(I))$, and $w \in \mathcal{W}^h$ is very similar to the one given in \cite{CS}, Propositions 4.1 and 4.2.  For the sake of brevity, we will refrain in reproducing it here.

We now prove that $\partial_3 \hat\vartheta=\partial_3 \vartheta$ in the sense of distributions,
where $\vartheta$ is the function entering in the definition of $u^1_\alpha$ in $\Omega\times D$.

For every $\varphi\in \mathcal{D}(\Omega\times Y)$ we have that
$$
\int_\Omega
\eps(\partial_1 u^\eps_2-\partial_2 u^\eps_1)\partial_3\varphi(x,\frac{x^\prime}\eps)
\,dx
=
2\eps\int_\Omega (Eu^\eps)_{23}\partial_1(\varphi(x,\frac{x^\prime}\eps))-
(Eu^\eps)_{13}\partial_2(\varphi(x,\frac{x^\prime}\eps))\,dx.
$$
Taking $\varphi(x,y)=\psi(x)\phi(y)$ with $\psi\in \mathcal{D}(\Omega)$ and 
$\phi\in \mathcal{D}(D)$ and evaluating the derivatives on the right-hand side of the identity above, by means of \eqref{tsvh20} and \eqref{tsvh6}, we may pass to the limit to find
$$
\int_{\Omega\times D}
(\partial_1 u_2+\partial^y_1 u_2^1-\partial_2 u_1-\partial^y_2 u_1^1)\phi\partial_3\psi
\,dxdy
=
2\int_{\Omega\times D}( K_{23}\partial_1^y\phi-
K_{13}\partial_2\phi)\psi
\,dxdy.
$$
Taking into account the structure of $u^1_\alpha$ and of $K_{\alpha3}$, the identity above
reduces to
$$
\int_{\Omega\times D}
2\vartheta\,\phi\,\partial_3\psi
\,dxdy
=
2\int_{\Omega\times D}-\partial_3\hat \vartheta\,\phi\,\psi
\,dxdy,
$$
that is: $\partial_3 \hat\vartheta=\partial_3 \vartheta$ in the sense of distributions.
Hence $\vartheta$ belongs to $L^2(\omega;H^1(I))$ and not just to $L^2(\Omega)$.
This implies that
$u^1_\alpha \in  L^2(\Omega; H^1_{\#}(Y))\cap L^2(\Omega; H^1_{m}(D))\cap L^2(\omega\times D;H^1(I))$. We can assume $u^1_\alpha$ to be equal to $0$ at $x_3=0$, since, as it follows from its definition \eqref{tsvh20}, $u^1_\alpha$ is defined up to a function of $x_3$. To remove this ambiguity we require $u^1_\alpha$ . This forces
$\vartheta=0$ at $x_3=0$, and as a consequence we deduce that $\vartheta=\hat \vartheta\in L^2(\omega;H^1_0(I))$. This, in turn implies that also $u^1_\alpha \in   L^2(\omega\times D;H^1_0(I))$.

Thus, to sum up, we have
\begin{equation}\label{K}
K=
\left(
\begin{array}{cc}
(E^y w)_{\alpha\beta} & \frac 12(y^R_\alpha \partial_3 \vartheta+\partial^y_\alpha v_3)\\
\frac 12(y^R_\alpha \partial_3 \vartheta+\partial^y_\alpha v_3) & \partial_3 u_3
\end{array}
\right).
\end{equation}

We now analyze \eqref{tsvh7}. The component $G_{33}$ is equal to zero since
$\dis \eps (E u^\eps)_{33} \chi_{M_\eps} \twoscale 0 $ by  \eqref{tsvh40}. 
For the component $G_{\alpha3}$ arising in \eqref{tsvh7},  we study the limit of each term in the sum 
$$ 2 \eps (E u^\eps)_{\alpha3} \chi_{M_\eps} = (\eps \partial_\alpha u_3^\eps + \eps  \partial_3 u_\alpha^\eps )\chi_{M_\eps}.$$
For the second term we use \eqref{tsvh100} (which implies that $\eps  u_\alpha^\eps \twoscale  u_\alpha$) to get
\begin{eqnarray}
\dis  \int_\Omega \eps  \partial_3 u_\alpha^\eps \phi(x) \psi(\frac{x'}{\eps}) \chi_{M_\eps}\ dx
& = &
- \int_\Omega \eps  u_\alpha^\eps \partial_3 \phi  \ \psi(\frac{x'}{\eps}) \chi_{M_\eps} \ dx \nonumber\\
& \to & 
- \int_\Omega \int_Y u_\alpha \partial_3 \phi(x)  \ \psi(y) \chi_{Y\setminus D} \ dx dy \label{secondterm} \\
& = &
  \int_\Omega \int_Y \partial_3 u_\alpha  \phi(x)  \ \psi(y) \chi_{Y\setminus D} \ dx dy, \nonumber
\end{eqnarray}
for all $\phi \in \mathcal{D}(\Omega)$ and $\psi \in C_{\#}({Y})$.
 By virtue of \eqref{tsvh50}, the two-scale limit of  $\eps \partial_\alpha u^\eps_3$ is equal to $ \partial^y_\alpha u_3$, thus
 $$
 2 \eps (E u^\eps)_{\alpha3} \chi_{M_\eps}\twoscale (\partial^y_\alpha u_3+\partial_3 u_\alpha)\chi_{Y\setminus D}=G_{\alpha3}\chi_{Y \setminus D}.
 $$
 To identify the components $G_{\alpha\beta}$ we use \eqref{tsvh20} to get
 \begin{eqnarray*}
2 \eps (E u^\eps)_{\alpha\beta} \chi_{M_\eps}& = &(\eps \partial_\alpha u^\eps_\beta +  \eps \partial_\beta u^\eps_\alpha) \chi_{M_\eps} \nonumber\\
&   \twoscale & (\partial_\alpha u_\beta + \partial_\alpha^y u_\beta^1 + 
  \partial_\beta u_\alpha + \partial_\beta^y u_\alpha^1)  \chi_{Y\setminus D} \\%\label{tsvh11}\\
  & = &
   2(Eu + E^y u^1)_{\alpha\beta} \chi_{Y\setminus D}=2 G_{\alpha\beta}\chi_{Y\setminus D}.\nonumber
 \end{eqnarray*} 
Hence
\begin{equation}\label{G}
G=
\left(
\begin{array}{cc}
(Eu + E^y u^1)_{\alpha\beta} & \frac 12(\partial^y_\alpha u_3+\partial_3 u_\alpha)\\
\frac 12(\partial^y_\alpha u_3+\partial_3 u_\alpha) & 0
\end{array}
\right),
\end{equation}
and the Lemma is proven.
\QED
 
 %%%%%%%%%%%%%%%%

\noindent
{\bf Proof of Theorem \ref{Mainresulthom}.} 
The convergences stated in the Theorem, with the strong two-scale convergences replaced by the weak two scale convergence, have been already proven in Lemma \ref{lem412}. Set 
$$
\sigma^\eps:=\Cbb Eu^\eps ,\quad 
\tau^\eps:=\eps\Cbb Eu^\eps 
$$
then, by \eqref{tensorhomogenization} and \eqref{bdh3} we deduce that
$$
\sigma^\eps\chi_{F_\eps} \twoscale \sigma\chi_D, \quad 
\tau^\eps \chi_{M_\eps} \twoscale \tau\chi_{Y\setminus D},
$$
for some $\sigma_{ij}, \tau_{ij}\in L^2(\Omega\times Y)$. For later use, we now prove that
\begin{equation}\label{sigmadiv}
\int_D\sigma_{\alpha\beta} (E^y\psi)_{\alpha\beta}+\sigma_{\alpha3} \partial_{\alpha}^y\psi_{3}\,dy=0
\end{equation}
almost everywhere in $\Omega$ for every $\psi\in H^1_{\#}(Y)$. 
Indeed, take $\varphi(x)=\phi(x)\psi(x^\prime/\eps)$, with $\phi\in \mathcal{D}(\Omega)$ and $\psi\in H^1_{\#}(Y)$, in \eqref{mainpbhomogenization} to find
$$
 \int_\Omega (\sigma^\eps\chi_{F_\eps}+\eps \tau^\eps\chi_{M_\eps})\cdot (\psi E\phi+\phi \frac 1\eps E^y\psi)\,dx =
\int_\Omega \eps f_\alpha(x, \frac{x'}{\eps})\cdot \varphi_\alpha(x) 
+ f_3(x, \frac{x'}{\eps})\cdot \varphi_3(x) \,dx.
$$
Multiplying the above identity by $\eps$ and by passing to the limit we deduce
$$
 \int_\Omega\int_Y \sigma(x,y)\cdot (\phi(x) E^y\psi(y))\chi_D\,dydx =0,
$$
which implies \eqref{sigmadiv}. In particular, 
\eqref{sigmadiv} implies that
\begin{equation}\label{sigmavinc}
\int_D \sigma_{\alpha\beta}\,dy=\int_D \sigma_{\alpha\beta} y_\gamma\,dy=0,
\end{equation}
as it can be checked by taking, in \eqref{sigmadiv}, $\psi_\alpha\in H^1_{\#}(Y)$ and equal to
$
\psi_\alpha=a_\gamma^\alpha y_\gamma+b_{\gamma\delta}^\alpha y_\gamma y_\delta
$
on $D$, where $a_\gamma^\alpha,b_{\gamma\delta}^\alpha\in\{0,1\}$.
Also, from \eqref{sigmadiv} we deduce
\begin{equation}\label{sigmavinc2}
\int_D \sigma_{\alpha3}\,dy=\int_D \sigma_{13} y_1\,dy=\int_D \sigma_{23} y_2\,dy=
\int_D \sigma_{13} y_2+\sigma_{23} y_1\,dy=0.
\end{equation}

 We are now in a position to exhibit the appropriate test function to pass to the limit in \eqref{mainpbhomogenization}. We set 
\begin{equation}\label{testhom}
\left\{ \begin{array}{ll} & \varphi_\alpha (x) = \frac{1}{\eps} \bar u_\alpha(x) + \bar u^1_\alpha(x, \frac{x'}{\eps})  + \eps \bar w_\alpha(x, \frac{x'}{\eps}), \\
& \varphi_3 (x) =  \bar u_3(x, \frac{x'}{\eps})+ \eps \bar v_3(x, \frac{x'}{\eps}),
\end{array} \right.
\end{equation}
where $\bar u_\alpha\in  \mathcal{D}(\Omega)$,  $\bar u^1_\alpha, \bar u_3 \in  \mathcal{D}(\Omega;H^1_{\#}({Y}))$, and $\bar w_\alpha, \bar v_3 \in  \mathcal{D}(\Omega)\times \mathcal{D}(\bar D)$.
Moreover, we assume that
\begin{eqnarray}
\bar u^1_\alpha(x,y)&=& \bar{\bar u}^1_\alpha(x)-y_\alpha^R\bar \vartheta(x)-y_\beta \partial_\beta \bar u_\alpha(x)\label{baru1},\\
\bar u_3(x,y)&=& \bar{\bar u}_3(x)-y_\beta \partial_3 \bar u_\beta(x)\label{baru3},
\end{eqnarray}
in $\Omega\times D$, where $\bar{\bar u}^1_\alpha, \bar{\bar u}_3, \bar \vartheta\in  \mathcal{D}(\Omega)$.
With this choice of test functions we find 
\begin{eqnarray*}
\eps (E\varphi)_{\alpha\beta}&=&(E\bar u+E^y\bar u^1)_{\alpha\beta}+\eps(E\bar u^1+E^y\bar w)_{\alpha\beta}+\eps^2(E^y\bar w)_{\alpha\beta},\\
2\eps (E\varphi)_{\alpha3}&=&(\partial_3\bar u_\alpha +\partial^y_\alpha\bar u_3)+\eps(\partial_3\bar u_\alpha^1+\partial_\alpha \bar u_3 +\partial^y_\alpha\bar v_3)+\eps^2(\partial_3\bar w_\alpha+\partial_\alpha \bar v_3),\\
\eps (E\varphi)_{33}&=&\eps\partial_3\bar u_3 +\eps^2\partial_3\bar v_3,
\end{eqnarray*}
and taking into account \eqref{baru1} and \eqref{baru3} we deduce that
\begin{eqnarray*}
(E\varphi)_{\alpha\beta}\chi_{F_\eps}&=&[(E\bar u^1+E^y\bar w)_{\alpha\beta}+\eps(E^y\bar w)_{\alpha\beta}]\chi_{F_\eps},\\
2(E\varphi)_{\alpha3}\chi_{F_\eps}&=&[(\partial_3\bar u_\alpha^1+\partial_\alpha \bar u_3 +\partial^y_\alpha\bar v_3)+\eps(\partial_3\bar w_\alpha+\partial_\alpha \bar v_3)]\chi_{F_\eps},\\
(E\varphi)_{33}\chi_{F_\eps}&=&[\partial_3\bar u_3 +\eps\partial_3\bar v_3]\chi_{F_\eps}.
\end{eqnarray*}
Problem \eqref{mainpbhomogenization} rewrites as
$$
\int_\Omega (\sigma^\eps\cdot E\varphi\chi_{F_\eps}+ \tau^\eps\cdot \eps E\varphi\chi_{M_\eps}\,dx =
\int_\Omega \eps f_\alpha(x, \frac{x'}{\eps})\cdot \varphi_\alpha(x) 
+ f_3(x, \frac{x'}{\eps})\cdot \varphi_3(x) \,dx,
$$
whose limit is
\begin{eqnarray}
&&\int_\Omega\int_D \sigma\cdot
\left(
\begin{array}{cc}
(E\bar u^1+E^y\bar w)_{\alpha\beta} & \frac 12(\partial_3\bar u_\alpha^1+\partial_\alpha \bar u_3 +\partial^y_\alpha\bar v_3)\\
\frac 12 (\partial_3\bar u_\alpha^1+\partial_\alpha \bar u_3 +\partial^y_\alpha\bar v_3) & \partial_3\bar u_3
\end{array}
\right)
\,dydx
\nonumber\\ 
& & \hspace{2cm}+
\int_\Omega\int_{Y\setminus D} \tau\cdot
\left(
\begin{array}{cc}
(E\bar u+E^y\bar u^1)_{\alpha\beta}  & \frac 12 (\partial_3\bar u_\alpha +\partial^y_\alpha\bar u_3)\\
\frac 12 (\partial_3\bar u_\alpha +\partial^y_\alpha\bar u_3) & 0
\end{array}
\right)
\,dydx\label{mainpbhom1}\\
& & \hspace{4cm}
=\int_\Omega \int_Yf_\alpha(x, y)\cdot \bar u_\alpha(x) 
+ f_3(x,y)\cdot \bar u_3(x,y) \,dydx.\nonumber
\end{eqnarray}
We claim that
\begin{equation}\label{mainpbhom2}
\int_\Omega\int_D \sigma_{\alpha\beta}
(E\bar u^1)_{\alpha\beta}\,dydx=0
\end{equation}
and that
\begin{equation}\label{mainpbhom3}
\int_\Omega\int_D\sigma_{\alpha 3}(\partial_3\bar u_\alpha^1+\partial_\alpha \bar u_3 +\partial^y_\alpha\bar v_3)\,dydx=\int_\Omega\int_D\sigma_{\alpha 3}(y_\alpha^R \partial_3\bar \vartheta+\partial^y_\alpha\bar v_3)\,dydx.
\end{equation}

Indeed, by means of \eqref{baru1} we compute
$$
(E\bar u^1)_{\alpha\beta}(x,y)= (E\bar{\bar u}^1)_{\alpha\beta}(x)-\frac 12 (y_\alpha^R\partial_\beta\bar \vartheta(x)+y_\beta^R\partial_\alpha\bar \vartheta(x))-y_\gamma \partial_\gamma (E{\bar u})_{\alpha\beta}(x)
$$
which, together with \eqref{sigmavinc}, implies \eqref{mainpbhom2}.
To prove \eqref{mainpbhom3}, note that, 
\begin{eqnarray*}
& &\partial_3\bar u_\alpha^1(x,y)+\partial_\alpha \bar u_3 (x,y)+\partial^y_\alpha\bar v_3(x,y)=\\
& &\hspace{1cm}[y_\alpha^R \partial_3\bar \vartheta(x)+\partial^y_\alpha\bar v_3(x,y)]
 +[\partial_3\bar{\bar u}^1_{\alpha}(x)+
\partial_\alpha\bar{\bar u}_{3}(x)]-y_\beta \partial_3 (E{\bar u})_{\alpha\beta}(x),
\end{eqnarray*}
where we used \eqref{baru1} and \eqref{baru3}.
Taking into account \eqref{sigmavinc2} we deduce that
\begin{eqnarray*}
& &\int_\Omega\int_D\sigma_{\alpha 3}[\partial_3\bar u_\alpha^1(x,y)+\partial_\alpha \bar u_3 (x,y)+\partial^y_\alpha\bar v_3(x,y)]\,dydx\\
& &\hspace{1cm}=\int_\Omega\int_D\sigma_{\alpha 3}[y_\alpha^R \partial_3\bar \vartheta(x)+\partial^y_\alpha\bar v_3(x,y)]-\sigma_{\alpha 3}y_\beta \partial_3 (E{\bar u})_{\alpha\beta}(x)\,dydx\\
& &\hspace{1cm}=\int_\Omega\int_D\sigma_{\alpha 3}[y_\alpha^R \partial_3\bar \vartheta(x)+\partial^y_\alpha\bar v_3(x,y)]-\sigma_{13}y_1 \partial_3 (E{\bar u})_{11}(x)\\
& &\hspace{2cm}-
\sigma_{23}y_2 \partial_3 (E{\bar u})_{22}(x)
-(\sigma_{13}y_2+\sigma_{23}y_1) \partial_3 (E{\bar u})_{12}(x)
\,dydx
\end{eqnarray*}
and using again \eqref{sigmavinc2} we obtain \eqref{mainpbhom3}.
Thanks to \eqref{mainpbhom2} and \eqref{mainpbhom3}, \eqref{mainpbhom1} reduces to
\begin{eqnarray}
&&\int_\Omega\int_D \sigma\cdot
\left(
\begin{array}{cc}
(E^y\bar w)_{\alpha\beta} & \frac 12(y_\alpha^R \partial_3\bar \vartheta+\partial^y_\alpha\bar v_3)\\
\frac 12 (y_\alpha^R \partial_3\bar \vartheta+\partial^y_\alpha\bar v_3) & \partial_3\bar u_3
\end{array}
\right)
\,dydx
\nonumber\\ 
& & \hspace{1cm}+
\int_\Omega\int_{Y\setminus D} \tau\cdot
\left(
\begin{array}{cc}
(E\bar u+E^y\bar u^1)_{\alpha\beta}  & \frac 12 (\partial_3\bar u_\alpha +\partial^y_\alpha\bar u_3)\\
\frac 12 (\partial_3\bar u_\alpha +\partial^y_\alpha\bar u_3) & 0
\end{array}
\right)
\,dydx\label{mainpbhom4}\\
& & \hspace{3cm}
=\int_\Omega \int_Yf_\alpha(x, y)\cdot \bar u_\alpha(x) 
+ f_3(x,y)\cdot \bar u_3(x,y) \,dydx,\nonumber
\end{eqnarray}
which holds for every $\bar u_\alpha\in  \mathcal{D}(\Omega)$,  $\bar u^1_\alpha, \bar u_3 \in  \mathcal{D}(\Omega;H^1_{\#}({Y}))$, $\bar w_\alpha, \bar v_3 \in  \mathcal{D}(\Omega)\times \mathcal{D}(\bar D)$, and $\bar \vartheta\in  \mathcal{D}(\Omega)$ that satisfy \eqref{baru1} and \eqref{baru3}.
By density \eqref{mainpbhom4} holds also for $(\bar u_i,\bar u^1_\alpha) \in \mathcal{U}^h$,  $(-y_2\bar\vartheta,y_1\bar\vartheta,v_3) \in \mathcal{V}^h$, and $w \in \mathcal{W}^h$.

To conclude the proof it suffices to show that
$$
\sigma\chi_{D}=\Cbb 
\left(
\begin{array}{cc}
(E^y w)_{\alpha\beta} & (E^{yx}v)_{\alpha3}\\
(E^{yx}v)_{\alpha3} & (Eu)_{33}
\end{array}
\right)\chi_{D}
$$
and
$$ 
\tau\chi_{Y\setminus D}=\Cbb 
\left(
\begin{array}{cc}
(Eu + E^y u^1)_{\alpha\beta} & (E^{yx}u)_{\alpha3}\\
(E^{yx}u)_{\alpha3} & 0
\end{array}
\right)
\chi_{Y\setminus D}.
$$
But, these identities follow essentially from the definition of two-scale convergence and Lemma \ref{lem412}. For instance, with $\varphi \in L^2(\Omega;  C_{\#}(Y))$ we have
$$
\int_\Omega \sigma^\eps_{ij}(x)\chi_{F_\eps}\varphi(x,\frac{x^\prime}\eps)\,dx= 
\int_\Omega \Cbb_{ijkl}(x,\frac{x^\prime}\eps)Eu^\eps_{kl}(x)\chi_{F_\eps}\varphi(x,\frac{x^\prime}\eps)\,dx
$$
and by passing to the limit we find
 $$
\int_\Omega\int_Y \sigma(x,y)\chi_{D}\varphi(x,y)\,dydx= 
\int_\Omega\int_Y \Cbb_{ijkl}(x,y)G(x,y)\chi_{D}\varphi(x,y)\,dydx
$$
where we used the fact that $\Cbb_{ijkl}\varphi\in L^2(\Omega;  C_{\#}(Y))$. This identity implies $\sigma\chi_{D}=\Cbb K\chi_{D}$.

The strong two-scale convergence of the rescaled strains can be proven as it was done in the rod problem.
\QED

%%%%%%%%%%%%%%%%

\noindent
{\bf Proof of Theorem \ref{nonlocaleffecthomhom}.} 
With
$$
z_\alpha:=u^1_\alpha+y_\beta \partial_3 u_\alpha-\vartheta y_\alpha^R,
\qquad
z_3:=u_3-(\mint_Du_3\,dy-y_\alpha\partial_3 u_\alpha)
$$
in $\Omega\times Y$,
it follows that the bilinear form of problem \eqref{Thelimitpbhom} defined on $\Omega\times (Y\setminus D)$ rewrites as
$$
\begin{array}{l}\displaystyle 
\int_\Omega\int_{Y\setminus D} 
\Cbb
\left(
\begin{array}{cc}
(Eu + E^y u^1)_{\alpha\beta} & (E^{yx}u)_{\alpha3}\\
(E^{yx}u)_{\alpha3} & 0
\end{array}
\right)
\\ \hspace{5cm}
\cdot
\left(
\begin{array}{cc}
(E\bar u+E^y\bar u^1)_{\alpha\beta}  &(E^{yx}\bar u)_{\alpha3}\\
(E^{yx}\bar u)_{\alpha3} & 0
\end{array}
\right)
\,dydx\\[10pt]
\displaystyle  
 =
\int_\Omega\int_{Y\setminus D} 
\Cbb
\left(
\begin{array}{cc}
(E^y z)_{\alpha\beta} & \frac 12 \partial^{y}_{\alpha}z_3\\
\frac 12 \partial^{y}_{\alpha}z_3 & 0
\end{array}
\right)
%\\ \hspace{5cm}
\cdot
\left(
\begin{array}{cc}
(E^y \bar z)_{\alpha\beta} & \frac 12 \partial^{y}_{\alpha}\bar z_3\\
\frac 12 \partial^{y}_{\alpha}\bar z_3 & 0
\end{array}
\right)
\,dydx.
\end{array}
$$
By means of the definitions of $u^1_\alpha$ and $u_3$ on $\Omega\times D$ we infer that
 $z_i=0$ in $\Omega\times D$. 
With this new notation we have that
$$
 U(x) := \mint_Y u_3\ dy =  \mint_Y z_3+(\mint_Du_3\,dy-y_\alpha\partial_3 u_\alpha)\, dy=\mint_Du_3\,dy+\mint_Y z_3\,dy.
$$
The proof of theorem is now very similar to the proof of the corresponding Theorem \ref{nonlocaleffect} in the rod problem. Since it can be achieved by following exactly the same steps, with the role of the variable $x'$ played by the variable $y$ and the spaces ${\mathcal Z},  {\mathcal Z}_0$ replaced by  ${\mathcal Z}^h, \mathcal{Z}^{h}_0 $, we refrain from writing all the details of the proof.
\QED

\noindent  {\Small  {\sc Acknowledgments.}
The first author gratefully acknowledges the hospitality and the support provided by the Institut de Math\'ematiques de Marseille
during the completion of this work.% in March 2015.
} 

%%%%%%%%%%%
\def\cprime{$'$} \def\cprime{$'$} \def\cprime{$'$}

\end{document}